\documentclass[reqno]{amsart}

\usepackage{amscd,amssymb,amsmath,latexsym,enumerate}
\usepackage{bbm}
\usepackage[mathscr]{euscript}
\usepackage{epsfig}
\usepackage{color}

\usepackage{a4wide}

\numberwithin{equation}{section}

\DeclareMathAlphabet{\mathpzc}{OT1}{pzc}{m}{it}

\newtheorem{theo}{Theorem}[section]
\newtheorem{defini}[theo]{Definition}
\newtheorem{proposi}[theo]{Proposition}
\newtheorem{lemma}[theo]{Lemma}
\newtheorem{coro}[theo]{Corollary}
\newtheorem{rem}[theo]{Remark}
\newtheorem{exam}[theo]{Example}

\newcommand{\CM}{{\mathbb C}}
\newcommand{\NM}{{\mathbb N}}

\newcommand{\RM}{{\mathbb R}}

\newcommand{\ZM}{{\mathbb Z}}

\newcommand{\Bb}{{\mathcal B}}
\newcommand{\Cc}{{\mathcal C}}
\newcommand{\Hh}{{\mathcal H}}

\newcommand{\Ll}{{\mathcal L}}

\newcommand{\Rr}{{\mathcal R}}

\newcommand{\tR}{{\widetilde{R}}}
\newcommand{\tr}{{\widetilde{r}}}

\newcommand{\htt}{{\widehat{t}}}

\newcommand{\as}{{\mathscr A}}

\newcommand{\js}{{\mathscr J}}
\newcommand{\ks}{{\mathscr K}}

\newcommand{\ns}{{\mathscr N}}

\newcommand{\AG}{{\mathfrak A}}
\newcommand{\CG}{{\mathfrak C}}
\newcommand{\OG}{{\mathfrak O}}

\newcommand{\zz}{\mathfrak{z}}

\newcommand{\Cs}{$C^{\ast}$-algebra }              
\newcommand{\CS}{$C^{\ast}$-algebra}               
\newcommand{\Css}{$C^{\ast}$-algebras }            
\newcommand{\CsS}{$C^{\ast}$-algebras}             
\newcommand{\Zd}{{{\mathbb Z}^{\mathit{d}}}}       
\newcommand{\dist}{\mbox{\rm dist}}                
\newcommand{\supp}{\mbox{\rm supp}}                
\newcommand{\Orb}{{\mathit Orb}}                   
\newcommand{\tra}{\mbox{\sc t}}                    
\newcommand{\tri}{\mbox{\tiny\sc t}}               

\newcommand{\Inv}{\js}								
\newcommand{\csp}{\as^{\Ll}}						
\newcommand{\Id}{\rm I}								
\newcommand{\qp}{\Psi}								
\newcommand{\Mqp}{\widehat{\qp}\!_z}				
\newcommand{\Mch}{\widehat{\chi}}  					
\newcommand{\dco}{d_{\csp}}        					
\newcommand{\Hdco}{\dco^H} 							
\newcommand{\Lip}{{\rm Lip}(\RM^d)} 				

\allowdisplaybreaks

\usepackage{color}

\begin{document}
\title{H\"older Continuity of the Spectra for Aperiodic Hamiltonians}
\author{Siegfried Beckus, Jean Bellissard, Horia Cornean}

\address{Institut f\"ur Mathematik\\
Universit\"at Potsdam\\
Potsdam, Germany}
\email{beckus@uni-potsdam.de}

\address{Westf\"alische Wilhelms-Universit\"at\\
Fachbereich 10 Mathematik und Informatik\\
Einsteinstrasse 62\\
48149, M\"unster, Germany}
\email{jeanbel@math.gatech.edu}

\address{Department of Mathematical Sciences\\
Aalborg University\\
Skjernvej 4A\\
DK-9220 Aalborg, Denmark 
Aalborg, Danemark}
\email{cornean@math.aau.dk}

\begin{abstract}
We study the spectral location of a strongly pattern equivariant Hamiltonians arising through configurations on a colored lattice. Roughly speaking, two configurations are "close to each other" if, up to a translation, they "almost coincide" on a large fixed ball. The larger this ball is, the more similar they are, and this induces a metric on the space of the corresponding dynamical systems. 
Our main result states that the map which sends a given configuration into the spectrum of its associated  Hamiltonian, is H\"older (even Lipschitz) continuous in the usual Hausdorff metric. Specifically, the spectral distance of two Hamiltonians is estimated by the distance of the corresponding dynamical systems. 
\end{abstract}


\maketitle


\section{Introduction}
\label{LipZd.Chap-Int}

\noindent This work is a follow-up on a series of papers concerning periodic approximations for Hamiltonians modelling aperiodic media \cite{BB16,Bec16,BBdN17,BBdN18,BP17}. Such Hamiltonians are bounded self-adjoint operators defined as effective models describing the behavior of conduction electrons in a solid. In the most common case the electrons are considered as independent. The first task in considering a material of this type is to compute the spectrum as a set and then the density of states (electronic properties). From a mathematical point of view, the next important task would be to determine the nature of the spectral measures, and then to investigate the transport properties.

\vspace{.1cm}

\noindent In the previous series of papers, systematic methods have been developed to compute the spectrum as a set through a sequence of periodic approximations \cite{BB16,Bec16,BBdN17,BBdN18}, as well as the density of states \cite{BP17}. In the present work, a special class of models is investigated for which the speed of convergence can be evaluated more accurately in terms of the distance of the associated dynamical systems. In order to do so, a metric is introduced on the space of possible atomic configurations and it will be proved that, at least for the class considered here, the spectrum, as a compact subset of the real line, is a H\"older (even Lipschitz) continuous function when the set of compact subsets of $\RM$ is endowed with the Hausdorff metric.

 \subsection{The class of models}
 \label{LipZd.ssect-fram}

The Euclidean distance on $\RM^d$ is denoted by $|x|:=\sqrt{\sum_{j=1}^d x_j^2}$ and the max norm on $\RM^d$ is given by $|x|_{\max}:=\max_{1\leq j\leq d} |x_j|$. 

\vspace{.1cm}

The class of systems considered here models a solid with atoms on a lattice $\Ll\subseteq\RM^d$.  Here a lattice is a subgroup that is isomorphic to $\ZM^d$, namely $\Ll=M\ZM^d$ where $M$ is a $d\times d$ invertible matrix with real coefficients. The properties of each atom are encoded by a letter in an alphabet $\as$, where $(\as,d_\as)$ is a compact metric space. The alphabet $\as$ may, for instance, encode the chemical species of atomic nuclei if $\as$ is finite. In addition, recent developments \cite{KP17} show that alphabets that are compact metric spaces (such as compact intervals) can be used in practice as a tool. 
The family of all atomic configurations is represented by the infinite Cartesian product $\csp:=\prod_{x\in\Ll}\as=\{\xi:\Ll\to\as\}$. Since $\as$ is compact, this space is compact. 

\vspace{.1cm}

Define the deformed cube $Q_r:=\{M x\in\RM^d\,:\, |x|_{\max}\leq r\}$ of side length $r>0$. The configuration space $\csp$ becomes a compact metric space if equipped with the metric
\begin{align}\label{eq:dco}
\dco (\xi,\eta):=
    \min\left\{\inf\left\{
      \frac{1}{r}\;:\; r\in(0,\infty)\,,\, d_\as\big(\xi(x),\eta(x)\big)\leq \frac{1}{r} \text{ for all } x\in Q_r\cap\Ll
        \right\}\,,\, 1\right\}\,,
\end{align}

\noindent c.f. Lemma~\ref{LipZd.Lem-Metr}. The infimum over the empty set equals $+\infty$ by convention. If $\as$ is finite we choose $d_\as$ to be the discrete metric. In this case, if $\dco (\xi,\eta)\leq\frac{1}{r}$ for some $r>1$, then $\xi|_{Q_{r}}=\eta|_{Q_{r}}$.

\vspace{.1cm}

The translation group $\Ll$ acts naturally on this space, as a full shift denoted by $\tra$, namely
$$\forall \xi\in \csp\,,\;
   \xi=\big(\xi(x)\big)_{x\in\Ll}\,,
    \hspace{1.5cm}
     \big(\tra^h\xi\big)(x):=\xi(x-h)\,,
      \hspace{1.5cm}
       h\in\Ll\,.
$$

 A quantum particle moving on the lattice $\Ll$, such as a valence electron, is modeled by:

\begin{itemize}
\item[(i)] a Hilbert space of states taken to be $\Hh:=\ell^2(\Ll)\otimes \CM^N$, where $N$ represents the number of internal degrees of freedom of the particle, and

\vspace{.1cm}

\item[(ii)] a bounded self-adjoint Hamiltonian $H_\xi$, acting on $\Hh$, where $\xi\in\csp$ represents the atomic configuration with which the particle interacts. 
\end{itemize}

\noindent For the sake of notation, $M_N(\CM)$ denotes the set of $N\times N$ matrices with complex coefficients equipped with the operator norm $\|M\|_{\rm op}:=\sup_{|x|\leq 1} |Mx|$ defined by the Euclidean norm. In addition, we use $\|M\|_{\max}:=\sup_{|x|_{\max}\leq 1} |Mx|_{\max}$ induced by the max-norm.

\vspace{.1cm}

\noindent Let $z\in\Ll$ be a translation vector. Denote by $U_z$ the unitary operator acting on $\Hh$, which is induced by the translation with $z\in\Ll$. Then translating the origin of coordinates in the lattice is equivalent to shifting the atomic configuration in the background, leading to the following {\em covariance condition}
$$U_z H_\xi U_z^{-1}=
   H_{\tri^z\xi}\,.
$$

\noindent In addition, it is natural to consider the situation in which the map $\xi\in\csp\mapsto H_\xi\in \Bb(\Hh)$ is strongly continuous \cite{Be86}. It is worth remarking that thanks to the covariance condition, the spectra of $H_\xi$ and $H_{\tri^z\xi}$ are identical \cite{Be86,BBdN17}. Hence the spectrum does not change along the {\em orbit} $\Orb(\xi):=\{\tra^z\xi\,:\,z\in\Ll\}$ of $\xi$. In addition, thanks to the strong continuity condition, if 
$\eta\in\overline{\Orb(\xi)}$, then $H_\eta$ has its spectrum contained in the spectrum of $H_\xi$, i.e., $\sigma(H_\eta)\subseteq\sigma(H_\xi)$.
These basic conditions are satisfied in particular if the Hamiltonian obeys two additional properties:

\begin{itemize}

  \item There exists a finite subset $\Rr\subseteq \Ll$, which is called {\em range}, and some   $N\times N$-matrix valued complex coefficients $t_{h,x}(\xi)$ such that
$$H_\xi\psi(x) =
   \sum_{h\in \Rr} t_{h,x}(\xi) \psi(x-h)\,.
$$

\noindent In this case we say that the Hamiltonian has {\em finite range $\Rr$}. The covariance condition imposes that $t_{h,x}(\xi)= t_{h,0}(\tra^{-x}\xi)$, hence from now on we will simply write $t_h$ instead of $t_{h,0}$. The strong continuity imposes that $t_{h}$ is continuous {\em w.r.t.} $\xi\in\Omega$. The self-adjointness of $H_\xi$ imposes that $h\in \Rr$ if and only if $-h\in \Rr$ and that $t_{-h}(\xi)= t_{h}(\tra^{h}\xi)^\ast$. Here $A^\ast$ denotes the adjoint of the matrix $A\in M_N(\CM)$. This gives
\begin{eqnarray}
\label{eq:Ham}
H_\xi\psi(x) &=&
  \sum_{h\in \Rr} t_h(\tra^{-x}\xi) \psi(x-h)\,,
\end{eqnarray}

 \noindent where $t_h:\csp\to M_N(\CM)$ is continuous, $\Rr=-\Rr$ is finite and $t_{-h}(\xi)= t_{h}(\tra^{h}\xi)^\ast$.
 
  \item If the alphabet $\as$ is finite, a function $t:\csp\to M_N(\CM)$ is called {\em strongly pattern equivariant} \cite{Ke03} (and correspondingly the family $H=(H_\xi)_{\xi\in\csp}$ is  called {\em strongly pattern equivariant Hamiltonian}) if there exists a radius $R_t\geq 1$ such that if two configurations $\xi$ and $\eta$ coincide on $Q_{R_t}$, then $t(\xi)=t(\eta)$. The radius $1\leq R_t<\infty$ is called the {\em radius of influence} of $t$. 

\vspace{.1cm}

\noindent In the present framework $(\as,d_\as)$ is allowed to be any compact metric space, thus requiring $t(\xi)=t(\eta)$ whenever $\xi$ and $\eta$ coincide on $Q_{R_t}$ is too restrictive. Instead, we say that  $t:\csp\to M_N(\CM)$ is {\em strongly pattern equivariant with $\beta$-H\"older continuous coefficients} if there exists a radius $R_t\geq 1$ and a constant $C_t\geq 1$ such that
\begin{align}\label{eq:PatEq}
\|t(\xi)-t(\eta)\|_{\rm op}\;
	\leq\; C_t \,\, \sup_{y\in Q_{R_t}\cap\Ll} d_\as\big(\xi(y),\eta(y)\big)^\beta.
\end{align}

\noindent Therefore, the family $H=(H_\xi)_{\xi\in\csp}$ will be called {\em strongly pattern equivariant with $\beta$-H\"older continuous coefficients} whenever each of its coefficients $t_h$ requires the previous regularity constraints. In this case,  the {\em radius of influence $R_H$ of $H$} is defined by 
\begin{align*}
	R_H:=\sup\{R_{t_h}\,:\, h\in\Rr\}\,,
\end{align*}

\noindent and we also introduce (see \eqref{eq:PatEq} for the definition of $C_t$) 
\begin{align*}
C_{hop}
	:=\sup_{h\in\Rr} C_{t_h}<\infty.
\end{align*}
\end{itemize}


It is worth noticing that for $\as\subseteq\RM$ finite, the usual discrete Schr\"odinger operator
\begin{align}\label{eq:SchrZd}
H_\xi\psi(x)=
   \sum_{h\in\Zd, \, |h|=1} \psi(x-h)+\xi(x)\psi(x)\,,
    \qquad 
     x\in\Zd\,,
\end{align}

\noindent  acting on $\ell^2(\Zd)$ is strongly pattern equivariant \cite{Be86} since the map $v:\as^{\Zd}\to\RM$ where $v(\xi)=\xi(0)$ induces a multiplicative real potential $v(\tra^{-x}\xi)=\xi(x)$. The first term is the discrete Laplace operator $\Delta$, which models the kinetic energy. The range of $\Delta$ is $\Rr_\Delta=\{\pm e_1,\ldots, \pm e_d\}$ where $\{e_1,\ldots, e_d\}$ is the standard basis of $\ZM^d$, while its radius of influence equals to $1$. Also, $v$ is "local" and only sees one point at the time, i.e. the origin. Thus $\Rr_v=\{0\}$ and its radius of influence is $R_v=1$. It follows that the range of $H_\xi$ is $\Rr=\Rr_\Delta\cup\Rr_v$ while its radius of influence $R_H=1$. 

 \subsection{A concise description of our main results}
 \label{LipZd.ssect-main}

\noindent We are mainly interested in the spectral properties of an operator family $H_\Xi:=(H_\xi)_{\xi\in\Xi}$. A subset $\Xi\subseteq\csp$ is called {\em invariant} if given any $\xi\in\Xi$ and $h\in\Ll$, then $\tra^h\xi\in\Xi$. An invariant, closed set $\Xi$ is called a {\em subshift} and the set of all subshifts is denoted by $\Inv$. Then $\Inv$  may be equipped with the Hausdorff metric $\Hdco$ induced by $\dco$, see \eqref{eq:HausDis}. Thus, $(\Inv,\Hdco)$ becomes a compact metric space, c.f. \cite[Proposition~4]{BBdN17}.

\vspace{.1cm}

For $\Xi\in\Inv$, the spectrum $\sigma(H_\Xi)$ of the operator family $H_\Xi:=(H_\xi)_{\xi\in\Xi}$ is defined by $\overline{\bigcup_{\xi\in\Xi}\sigma(H_\xi)}$. Since 
$\|H_\xi\|$ is uniformly bounded in $\xi\in\Xi$ and $H_\xi$ is self-adjoint, $\sigma(H_\Xi)$ is a compact subset of $\RM$. Then the distance between two spectra $\sigma(H_\Xi)$ and $\sigma(H_\Theta)$ for $\Xi,\Theta\in\Inv$ is defined by the Hausdorff metric $d_H$ on the set $\ks(\RM)$ of compact subsets of $\RM$ \cite{Ha62}, see Section~\ref{LipZd.ssect-ham} for details.

\vspace{.1cm}

\noindent In recent works \cite{Bec16,BBdN17} it has been shown that the map 
$$\Sigma_H:\Inv\to\ks(\RM)
	\,,\;
	\Xi\mapsto\sigma(H_\Xi)\,,
$$

\noindent is continuous in the corresponding Hausdorff topologies where $H=(H_\xi)_{\xi\in\Xi}$ is a Hamiltonian as defined in \eqref{eq:Ham}. We note that the continuity result holds under much more general conditions than we consider here and its proof relies on constructing a continuous field of {\CsS} using the uniform continuity of the operator coefficients.

\vspace{.1cm}

\noindent Our current paper deals with the question regarding the connection between the regularity of the coefficients (see \eqref{eq:PatEq}) and the regularity of the map $\Sigma_H$.

\vspace{.1cm}

\noindent The decay of the off-diagonal influences additionally the regularity of the spectra. This is measured by the Schur-$\beta$ norm $\|\cdot\|_\beta$ for $0<\beta\leq 1$. Specifically, this norm is defined for a finite range Hamiltonian $H=(H_\xi)_{\xi\in\csp}$ given in \eqref{eq:Ham} by
\begin{align}\label{eq:Schur}
\|H\|_\beta
	:= \sum_{h\in\Rr} \|t_h\|_\infty (1+|h|^2)^{\frac{\beta}{2}}\,,
	\quad\text{ where }\quad
	\|t_h\|_\infty:=\sup_{\xi\in \csp}\|t_h(\xi)\|_{\rm op}.
\end{align}

\noindent For simplification, let us first present our main result in the special case where $\as$ is finite. It asserts that $\Sigma:\Inv\to\ks(\RM),\, \Xi\mapsto\sigma(H_\Xi)$ is Lipschitz continuous for every strongly pattern equivariant Hamiltonian $H$.

\begin{proposi}
\label{LipZd.Prop-AlpFin}
Let $\as$ be finite and consider a strongly pattern equivariant Hamiltonian $H=(H_\xi)_{\xi\in\csp}$ with finite range $\Rr$ and strongly pattern equivariant coefficients $t_h,\, h\in\Rr$ as defined in \eqref{eq:Ham}. Then there exists a constant $C_{d,\Ll}$  such that $\Sigma_H$ is Lipschitz continuous:
$$d_H\big(\sigma(H_\Xi),\sigma(H_\Theta)\big)\leq 
   C_{d,\Ll}\; C_{hop}\; R_H\; \|H\|_1\; \Hdco\big(\Xi,\Theta\big)
   \,,\qquad \Xi,\Theta\in\Inv\,.
$$
\end{proposi}

\noindent The latter result is an immediate consequence of the following main theorem where $(\as,d_\as)$ is a compact metric space only.

\begin{theo}
\label{LipZd.th-main}
Let $H=(H_\xi)_{\xi\in\csp}$ be a finite range strongly pattern equivariant Hamiltonian with $\beta$-H\"older continuous coefficients for some $0<\beta\leq 1$. Then there exists a constant $C_{d,\Ll}$ only depending on the dimension $d$ and the lattice $\Ll$ such that $\Sigma_H$ is $\beta$-H\"older continuous:
$$d_H\big(\sigma(H_\Xi),\sigma(H_\Theta)\big)\leq 
   C_{d,\Ll}\; C_{hop} \; \|H\|_\beta\; R_H^\beta\; \Hdco\big(\Xi,\Theta\big)^\beta,
   \qquad \Xi,\Theta\in\Inv\,.
$$
\end{theo}

\vspace{.1cm}

\noindent {\bf Proof of Proposition~\ref{LipZd.Prop-AlpFin}:} If $\as$ is finite, every strongly pattern equivariant function $t_h:\csp\to\CM$ is Lipschitz continuous ($\beta=1$) with Lipschitz constant $C_{t_h}$ and some $R_{t_h}\geq 1$. Thus, $C_{hop}=\max\{C_{t_h}\,:\, h\in\Rr\}$ and $R_H=\max\{R_{t_h}\,:\, h\in\Rr\}$ are finite. Hence, the statement follows from Theorem~\ref{LipZd.th-main}.
\hfill$\Box$

\vspace{.2cm}

\begin{exam}
\label{LipZd.Ex-AlpFin}
{\em Consider the Schr\"odinger operator defined in \eqref{eq:SchrZd} with nearest neighbor interaction on $\Ll=\Zd$. Clearly, Proposition~\ref{LipZd.Prop-AlpFin} implies the Lipschitz continuity of the spectra with $R_H=1$, see discussion right after \eqref{eq:SchrZd}.
}
\hfill $\Box$
\end{exam}

The previous spectral estimates can be extended to Hamiltonians with infinite range. Specifically, let $H=(H_\xi)_{\xi\in\csp}$ be given by \eqref{eq:Ham} while the range $\Rr$ is infinite. For $\beta>0$, we call $H$ a {\em strongly pattern equivariant Hamiltonian with $\beta$-H\"older continuous coefficients and infinite range} if $\|H\|_\beta$ and $C_{hop}$ are both finite. Let $H|_r$ be the restriction of $H$ to the range $\Rr\cap Q_r$, see \eqref{eq:RestHam} for details. If the radius of influence $R_{H|_r}$ satisfies 
\begin{align*}
R_{H|_r} \leq C_H\; r\,,\qquad
	r\geq 1\,
\end{align*}

for some constant $1\leq C_H<\infty$ independent of $r$, we say that $H$ admits a {\em radius of influence with linear growth} in $r$. Then we have the following result:

\begin{theo}
\label{LipZd.th-infirange}
Let $H=(H_\xi)_{\xi\in\csp}$ be a strongly pattern equivariant Hamiltonian (possible infinite range) with $\beta$-H\"older continuous coefficients for $0<\beta\leq 1$. If the radius of influence of $H$ has a linear growth with constant $1\leq C_H<\infty$, then
$$d_H\big(\sigma(H_\Xi),\sigma(H_\Theta)\big) \;
	\leq \; 2\; C_{d,\Ll} \; \|H\|_\beta \; (C_H^\beta+C_{hop}) \; \Hdco\big(\Xi,\Theta\big)^\beta
	\,,\qquad \Xi,\Theta\in\Inv\,,
$$

where $C_{d,\Ll}>0$ is the same constant as in Theorem~\ref{LipZd.th-main}.
\end{theo}

\vspace{.1cm}

\noindent For $\xi\in\csp$, consider $\Xi_\xi:=\overline{\Orb(\xi)}\in\Inv$. Such subshifts are called {\em topological transitive}. Note that the set of topological transitive subshifts $\{\Xi_\xi\,:\, \xi\in\csp\}\subseteq\Inv$ is not closed in the Hausdorff topology \cite[Example~2]{BBdN18}.

\begin{coro}
\label{LipZd.Cor-toptrans}
Let $H=(H_\xi)_{\xi\in\csp}$ be a strongly pattern equivariant Hamiltonian (possibly with infinite range) having $\beta$-H\"older continuous coefficients where $0<\beta\leq 1$. If the radius of influence of $H$ has linear growth with constant $1\leq C_H<\infty$, then 
$$d_H\big(\sigma(H_\xi),\sigma(H_\eta)\big) \;
	\leq \; 2\; C_{d,\Ll} \; \|H\|_\beta \; (C_H^\beta+C_{hop}) \; \Hdco\big(\Xi_\xi,\Xi_\eta\big)^\beta
	\,,\qquad \xi,\eta\in\csp\,,
$$

where $C_{d,\Ll}>0$ is the same constant as in Theorem~\ref{LipZd.th-main}.
\end{coro}

\vspace{.1cm}

\noindent {\bf Proof: } Since $\xi\mapsto H_\xi$ is continuous in the strong operator topology and the operator family is equivariant, the inclusion $\sigma(H_\eta)\subseteq\sigma(H_\xi)$ follows for all $\eta\in\Xi_\xi$. Hence, $\sigma(H_{\Xi_\xi})=\sigma(H_\xi)$ is derived and Theorem~\ref{LipZd.th-infirange} finishes the proof.
\hfill$\Box$

\vspace{.2cm}

\noindent In particular, if $\Xi_\xi=\Xi_\eta$ holds for some $\xi,\eta\in\csp$, then their spectra $\sigma(H_\xi)$ and $\sigma(H_\eta)$ coincide for any strongly pattern equivariant Hamiltonian. Note that $\Xi_\xi=\Xi_\eta$ defines an equivalence relation on $\csp$, which we call {\em orbit closure equivalence}.

\vspace{.1cm}

\noindent Clearly, Lipschitz/H\"older continuity of $\Sigma_H$ implies the continuity of this map. This provides a different proof (in a special case) of the much more general statement \cite[Theorem~2]{BBdN17}, namely

\begin{coro}
\label{LipZd.Cor-ContSpec}
Let $H=(H_\xi)_{\xi\in\csp}\in\CG$ be a strongly pattern equivariant Hamiltonian (possible infinite range) with $\beta$-H\"older continuous coefficients. Then the map
$$\Sigma:\Inv\to\ks(\RM)\,,\qquad
	\Xi\mapsto\sigma(H_\Xi)\,,
$$

is continuous in the corresponding Hausdorff topologies.
\end{coro}

\noindent {\bf Proof: } If $H$ has finite range, the statement follows from Theorem~\ref{LipZd.th-main}. If $\Rr$ is infinite, the statement follows by the fact that $H$ can be approximated in norm by $H_r$ with finite range.
\hfill$\Box$

 \subsection{Comments and Method}
 \label{LipZd.ssect-cometh} Continuity of spectral gaps has been proven in some cases in the past. The problem occurred, in particular, with the dependence of the spectrum as a function of an external, uniform magnetic field \cite{Si82, CP12, CP15}. It also occurred for the Schr\"odinger operator on the line $\RM$, with almost periodic potentials in which the frequency module is varying \cite{El82}. Lipschitz continuity of the spectral gaps was studied more thoroughly in \cite{Be94} for a broad class of models for $2D$-crystal electrons in a uniform magnetic field. This was extended in \cite{MP02,MP05,MPR05} as well. The problem of continuity {\em w.r.t.} changing the underlying atomic configurations was only recently systematically investigated  in connection with the need to compute the spectrum of a Hamiltonian describing the electron motion in an aperiodic environment \cite{BB16,Bec16,BBdN17,BBdN18}. In addition, the convergence of the density of states measures is studied in \cite{BP17}. For it has been remarked for a long time by physicists that periodic approximations provided the most efficient numerical method to achieve the result \cite{Ho76,OK85,SB91}. A folklore theorem was also showing that the accuracy of such methods was exponentially fast in the period of the approximation. This was proved in some way in \cite{Pr12,Pr17}. The existence and construction of periodic approximations in dimension one was recently characterized in \cite{BBdN18}. In several of the studies mentioned above, the concept of continuous field of \Css was explicitly used. One of the advantages of this concept is that the spectrum of a self-adjoint continuous section of such a field is always continuous \cite{Ka51,DD63,Di69,BB16}. However the question of whether a continuous field of \Css is Lipschitz continuous or even differentiable, has not been addressed in a systematic way so far. Therefore, to the best of their knowledge, the authors believe that the Lipschitz/H\"older continuity dependence of the spectrum in terms of the underlying atomic configuration, described in the present article, is new. 

\vspace{.1cm}

\noindent A bit of explanation for the method used here is in order, a method which employs and adapts a numbers of ideas from \cite{CP12, CP15} as well as \cite{BBdN17}. Consider the simplest case, the Schrödinger operator defined in \eqref{eq:SchrZd} by the discrete Laplacian plus a potential on $\Zd$. 
It is well-known, that if two different potentials are close in the uniform norm (sup norm), then their spectra is close as well in the Hausdorff metric. This follows as then their operator norm difference is small as well and all considered operators are self-adjoint. In the present paper, these potentials are described via different configurations in $\as^{\Zd}$. Clearly, the assumption that two potentials are close in the uniform norm is too restrictive in general. For instance, if the alphabet $\as$ is finite it yields that both potentials are equal if they are close enough in the uniform norm. This is in particular problematic if one seeks to approximate non-periodic configurations via periodic ones. To overcome this difficulty, the first powerful concept is coming from using the Hausdorff metric $\Hdco$ on the set of subshifts \cite{BBdN17}. 
Two configurations $\xi, \eta\in\Omega$ generate "close" subshifts if there exists a large radius $r$ such that they share almost the same local pattern of size $r$ \cite{Bec16,BBdN18}. Specifically, locally on large areas the two configurations are close modulo translation.
The second method \cite{CP12, CP15} consists in showing that there is $\delta >0$ and $C>0$ such that if a complex number $\zz$ satisfies $\dist(\zz,\sigma(H_\eta))>C\Hdco(\Xi_\xi,\Xi_\eta)^{\delta}$, then $\zz$ belongs to the resolvent set of $H_\xi$ as well. This is obtained through localizing the resolvent operators with the help of a Lipschitz-partition of unity. This trick implements the fact that the potentials are only locally close but not uniformly as discussed before. Then for each local region selected this way, in which $H_\xi$ sees a local pattern, shared with $\eta$, there is a translation at finite distance bringing the same pattern for $\eta$ in this local region. Hence up to a translation depending on the region of localization, $H_\eta$ and $H_\xi$ are close in this region, thanks to the strongly pattern equivariance condition. This closeness is translated into a proof of the existence of $(\zz-H_\xi)^{-1}$, see Proposition~\ref{LipZd.Prop-EsRes}. 

\vspace{.1cm}

\noindent This result is actually coming as a surprise. Indeed, in view of \cite{BB16}, the best that could be expected is, for a Hamiltonian with Lipschitz continuous coefficients, to be $p2$-Lipschitz. In such a case, the Hausdorff distance between spectra should only be H\"older continuous of exponent $1/2$. As discussed in \cite{BB16}, such a loss of regularity is usually due to gap closing, as observed, for instance, in the Harper model \cite{Ho76}. In the present case, the pattern equivariance condition is actually a strong constraint on the system. Theorem~\ref{LipZd.th-main} suggest then that, even if there are gaps closing in some limit, it cannot occur at a slower rate than the one imposed by this Lipschitz continuity. One model has been numerically investigated in the literature, which is called the Kohmoto model \cite{OK85}. It represents a paradigm for one-dimensional quasicrystals. It contains a real parameter $\alpha$ that labels the slope of the line implementing the physical space in a cut-and-project scheme. It will be shown in a forthcoming paper \cite{BB19}, that the present estimate applies to this parameter, in that $\alpha$ defines a specific subshift $\Xi_\alpha$ to which Theorem~\ref{LipZd.th-main} applies. The combinatoric distance implements, on the set of slopes, a topology making the real line completely disconnected, but compatible with the encoding by a continuous fraction expansion. It is actually difficult from looking at the numerical results in \cite{OK85} to see this Lipschitz continuity, because the map $\alpha\in\RM \to \Xi_\alpha\in\js$ is actually discontinuous if $\RM$ is endowed with its usual Euclidean topology \cite{BIT91}.

\vspace{.1cm}

\noindent In light of \cite{BBdN17,BP17}, it is natural to ask for extensions for Hamiltonians on Delone sets. This setting includes important geometric examples such as the Penrose tilings or the octagonal lattice, which are not included here. The main technical difficulty is that the corresponding operators for different Delone sets act on different Hilbert spaces that are not isomorphic, in general. This problem is studied in a forthcoming work. Furthermore, the continuity result of the map $\Sigma_H$ in \cite{BBdN17} requires amenability. Since some arguments provided here extend to more general groups, it is interesting to ask if there is a connection between the amenability and the existence of a Lipschitz partition of unity in the group.

 \subsection{Organization of the paper}
 \label{LipZd.ssect-orga}

More details about the configuration space and the Hausdorff topology are provided in Section~\ref{LipZd.ssect-conf}. We discuss several properties of the considered class of Hamiltonians in Section~\ref{LipZd.ssect-ham}. Section~\ref{LipZd.Chap-tech} is devoted to the concept of Lipschitz-partitions of unity and superoperators that are crucial ingredients in the proof of the main theorem. In Section~\ref{LipZd.Chap-proof}, the proofs of the main results are given. 

\vspace{.1cm}

\subsection*{Acknowledgements} This research was supported through the program ``Research in Pairs" by the Mathematisches Forschungsinstitut Oberwolfach in 2018. This research has been supported by grant 8021--00084B \emph{Mathematical Analysis of Effective Models and Critical Phenomena in Quantum Transport} from The Danish Council for Independent Research \textbar\ Natural Sciences.


\section{Framework}
\label{LipZd.Chap-Fra}

\noindent In this section, an introduction is provided of the used concepts and notations of the configuration space $\csp$ and the Hausdorff metric $\Hdco$ on the set $\Inv$ of closed, invariant subsets of $\csp$. In the second part, a detailed elaboration of the studied Hamiltonians is presented. In addition, some auxiliary statements are provided that are used in the proof of the main theorem.

 \subsection{The configuration space and local patterns}
 \label{LipZd.ssect-conf}
 
Throughout this work $(\as,d_\as)$ is a compact metric space and the configuration space $\as^\Ll=\{\xi:\Ll\to\as\}$ is endowed with the product topology.

\begin{lemma}
\label{LipZd.Lem-Metr}
Let $(\as,d_\as)$ be a compact metric space. Then $(\csp,\dco)$ is a compact metric space inducing the product topology on $\csp$ where $\dco$ is defined in \eqref{eq:dco}
\end{lemma}

\noindent {\bf Proof: } 
Since $\as$ is compact, $\as^\Ll$ is compact as well. Thus it suffices to show that $\dco$ defines a metric as it is immediate to see that $\dco$ generates the product topology on $\csp$. Since $d_\as$ is a metric, it is straightforward to show $\dco(\xi,\eta)=\dco(\eta,\xi)$ and that $\dco(\xi,\eta)=0$ implies $\xi=\eta$. In order to prove the triangle inequality, let $\varepsilon>0$. Then there is an $r_\varepsilon>0$ and $\tr_\varepsilon>0$ satisfying 
$$
\dco(\xi,\eta) \leq \frac{1}{r_\varepsilon} \leq \dco(\xi,\eta) + \varepsilon
	\,,\qquad
\dco(\eta,\zeta) \leq \frac{1}{\tr_\varepsilon} \leq \dco(\eta,\zeta) + \varepsilon\,,
$$

and
$$
d_\as\big(\xi(x),\eta(x)\big)\leq \frac{1}{r_\varepsilon} 
	\,\text{ for all }\, x\in Q_{r_\varepsilon}\cap \Ll
	\,,\qquad
d_\as\big(\eta(y),\zeta(y)\big)\leq \frac{1}{\tr_\varepsilon}
	\,\text{ for all }\, y\in Q_{\tr_\varepsilon}\cap \Ll\,.
$$

Define $1/r := 1/r_\varepsilon + 1/\tr_\varepsilon$. Thus, $r<r_\varepsilon$ and $r<\tr_\varepsilon$ holds implying $Q_r\subseteq Q_{r_\varepsilon}$ and $Q_r\subseteq Q_{\tr_\varepsilon}$. Then the triangle inequality for $d_\as$ gives 
$$d_\as\big(\xi(x),\zeta(x)\big)
	\leq d_\as\big(\xi(x),\eta(x)\big) + d_\as\big(\eta(x),\zeta(x)\big)
	\leq \frac{1}{r_\varepsilon} + \frac{1}{\tr_\varepsilon}\,,\qquad
	x \in Q_{r}\cap \Ll\,.
$$

Hence, 
$$\dco\big(\xi,\zeta\big)
	\leq \frac{1}{r_\varepsilon} + \frac{1}{\tr_\varepsilon}
	\leq \dco\big(\xi,\eta\big) + \dco\big(\eta,\zeta\big) + 2\varepsilon
$$

follows for $\varepsilon>0$ arbitrary proving the triangle inequality for $\dco$.
\hfill$\Box$

\vspace{.2cm}

Since $\csp$ is compact, $\js$ is contained in $\ks\big(\csp\big)$ which denotes the set of compact subsets of $\csp$. The set $\ks\big(\csp\big)$ gets a compact metrizable space if equipped with the Hausdorff topology \cite{CV77}, also called Chabauty-Fell topology \cite{Ch50,Fe62}. More precisely, the Hausdorff metric $\Hdco:\ks\big(\csp\big)\times\ks\big(\csp\big)\to[0,\infty)$ induced by $\dco$ is defined as follows
\begin{align}\label{eq:HausDis}
\Hdco(\Xi_1,\Xi_2) := 
	\max\left\{
		\sup_{\xi_1\in\Xi_1} 
			\inf_{\xi_2\in\Xi_2} \dco(\xi_1,\xi_2)
		,\;
		\sup_{\xi_2\in\Xi_2} 
			\inf_{\xi_1\in\Xi_1} \dco(\xi_1,\xi_2)
	\right\}
\end{align}

\noindent and $\big(\ks(\csp),\Hdco\big)$ is a compact metric space. Since $\Inv\subseteq\ks\big(\csp\big)$ is a closed subset \cite{Bec16,BBdN17}, $(\Inv,\Hdco)$ is a compact metric space. If $\as$ is finite, the topology on $\Inv$ can also be described by its local pattern topology defined in \cite[Section~3.3]{Bec16}. It asserts that the convergence of subshifts in the Hausdorff topology is equivalent to the convergence of the local patterns \cite[Theorem~3.3.22]{Bec16}. This fact is implicitly used in this work (as the operators are localized on patches by the Lipschitz-partition of unity), which is implemented via the following lemma.

\begin{lemma}
\label{LipZd.Lem-BaEq}
The following assertions hold.
\begin{itemize}
\item[(a)] Let $\Xi,\Theta\in\Inv$. For every $\xi\in\Xi$ there exists an $\eta:=\eta(\xi)\in\Theta$ such that
$$\dco\big(\xi,\eta(\xi)\big)
	\;\leq\; \Hdco(\Xi,\Theta)\,.
$$
\item[(b)] If $\xi,\eta\in\csp$ satisfy $\dco(\xi,\eta)\leq\frac{1}{r}$ with $r>2$, then 
$$
d_\as\big(\xi(x),\eta(x)\big) \leq \frac{1}{r-1}\,,\qquad x\in Q_{r-1}\cap\Ll\,.
$$
\end{itemize}
\end{lemma}

\noindent {\bf Proof: }
(a) Let $\xi\in\Xi$ be arbitrary. Because $\Theta$ is compact, there exists $\eta(\xi)\in \Theta$ such that  $\dco(\xi,\eta(\xi))=\inf_{\eta\in\Theta} \dco(\xi,\eta)$. Then \eqref{eq:HausDis} implies $\dco(\xi,\eta(\xi))\leq \Hdco(\Xi,\Theta)$.

\vspace{.1cm}

(b) Since $r-1>1$ and $\dco(\xi,\eta)<\frac{1}{r-1}$ are strict, the estimate follows immediately by the definition of the infimum in \eqref{eq:dco}.
\hfill $\Box$

\vspace{.2cm}

\begin{rem}
\label{LipZd.Rem-BaEq}
{\em
(i) Clearly, the role of $\Xi$ and $\Theta$ can be interchanged. Specifically, for each $\eta\in\Theta$, there is an $\xi(\eta)\in\Xi$ such that $\dco\big(\xi(\eta),\eta\big) \,\leq \, \Hdco(\Xi,\Theta)$ .

\vspace{.1cm}

\noindent (ii) If $\as$ is finite, $d_\as$ is the discrete metric and Lemma~\ref{LipZd.Lem-BaEq} can be reformulated as follows:
Let $\Xi,\Theta\in\Inv$ be such that $r:=\Hdco(\Xi,\Theta)^{-1}>1$. Then for every $\xi\in\Xi$ there is an $\eta:=\eta(\xi)\in\Theta$ such that $\xi|_{Q_r}=\eta|_{Q_r}$. The reader is referred for a purely topological discussion of this observation in \cite{Bec16,BBdN18}. 
}
\hfill $\Box$
\end{rem}

 \subsection{Hamiltonians}
 \label{LipZd.ssect-ham}

\noindent As described in the introduction, the Hamiltonians are self-adjoint, bounded operators on the Hilbert space $\Hh:=\ell^2(\Ll)\otimes\CM^N$. Denote by $\Bb(\Hh)$ the \Cs of bounded, linear operators on the Hilbert space $\Hh$ equipped with the operator norm $\|\cdot\|$. The group $\Ll$ is represented by its left regular representation defined by
$$U_z\varphi(x) := \varphi(x-z)
	\,,\qquad \varphi\in\Hh\,,\, x,z\in\Ll\,.
$$

\noindent Let $\Rr\subseteq\Ll$ be a finite subset and $t_h:\csp\to M_N(\CM)$ be continuous for $h\in\Rr$. Define the hopping function $t_{h,\xi}:\Ll\to M_N(\CM)$ by $t_{h,\xi}(x):=t_h(\tra^{-x}\xi)$ for $\xi\in\csp$. With this at hand, the Hamiltonian $H_\xi:\Hh\to\Hh$ defined in Equation~\eqref{eq:Ham} is represented by
$$H_\xi:=\sum_{h\in\Rr} \htt_{h,\xi} \, U_h
$$

\noindent where $\widehat{f}:\Hh\to\Hh$ denotes the multiplication operator $(\widehat{f}\varphi)(x):=f(x)\varphi(x)$ by the function $f:\Ll\to M_N(\CM)$. The operator $H_\xi$ is linear and uniformly bounded in $\xi\in\csp$, namely
$$\|H\|:=\sup_{\xi\in\csp}\|H_\xi\| 
	\;\leq\; \sum_{h\in\Rr} \|t_h\|_\infty 
	\;<\; \infty\,.
$$

\noindent Since $t_h:\csp\to M_N(\CM)$ and $\tra:\csp\to\csp$ are continuous, $\xi\mapsto H_\xi$ is continuous with respect to the strong operator topology. Furthermore, an elementary computation leads to the covariance condition $U_z H_\xi U_z^{-1}=H_{\tra^z\xi}$ for $z\in\Ll$ and $\xi\in\csp$.

\vspace{.1cm}

\noindent Let $\AG$ be the $\ast$-algebra generated by the operator families $(H_\xi)_{\xi\in\csp}$ defined in Equation~\eqref{eq:Ham}. The involutive and algebraic structure on $\AG$ is defined pointwise in $\xi\in\csp$. Equipped with the norm $\|H\|:=\sup_{\xi\in\csp}\|H_\xi\|$, $\AG$ is a normed $\ast$-algebra and its completion $\CG$ is a \CS. Then $\CG$ is a sub-\Cs of the reduced groupoid \Cs defined by the transformation groupoid $\csp\rtimes_{\tra}\Ll$, see e.g. \cite{Bec16,BBdN17}. If $\Xi\in\Inv$ is a subshift, a \Cs $\CG(\Xi)$ is similarly defined. Specifically, it is the closure of the $\ast$-algebra of operators $H_\Xi:=(H_\xi)_{\xi\in\Xi}$. This \Cs is again a sub-\Cs of the reduced groupoid \Cs defined by $\Xi\rtimes_{\tra}\Ll$. 

\vspace{.1cm}

\noindent The main focus in this work is on the study of self-adjoint $H\in\CG$ with strongly pattern equivariant hopping functions. Before providing the precise definition, the following auxiliary statement is proven.

\begin{lemma}
\label{LipZd.Lem-BasHam}
Let $h\in\Ll$ and $t:\csp\to M_N(\CM)$ be continuous.
\begin{itemize}
\item[(a)] Then $\|t\|_\infty=\|t\circ\tra^{-h}\|_\infty$ and $U_{-h}\, \htt_\xi = \htt_{\tra^{-h}\xi}\, U_{-h}$ hold for every $\xi\in\csp$.
\item[(b)] If, additionally, $t:\csp\to M_N(\CM)$ is a strongly pattern equivariant function with $\beta$-H\"older continuous coefficients with constant $C_t$ and radius of influence $R_t$, then $t\circ\tra^{-h}:\csp\to\CM$ is strongly pattern equivariant with $\beta$-H\"older continuous coefficients with radius of influence $R_t+|M^{-1}h|_{\max}$ and constant $C_{t\circ\tra^{-h}}=C_t$.
\end{itemize}
\end{lemma}

\noindent {\bf Proof: } (a) The equality $\|t\|_\infty=\|t\circ\tra^{-h}\|_\infty$ follows immediately by definition of the uniform norm, c.f. Equation~\eqref{eq:Schur}. Let $\xi\in\csp$. Then the identities
$$\big(U_{-h} \, \htt_\xi \, \varphi\big)(x) \;
	= \; t\big(\tra^{-x-h}\xi\big)\cdot \varphi(x+h)\;
	= \; \left(\htt_{\tra^{-h}\xi} \, U_{-h}\,\varphi\right)(x)
$$

hold for every $\varphi\in\Hh$ and $x\in\Ll$.

\vspace{.1cm}

\noindent (b) Since the inclusion
$$Q_{R_t}+h \;
	=\; M\big\{ \tilde{y}+M^{-1}h\in\RM^d \,:\, |\tilde{y}|_{\max}\leq R_t \big\} \;
	\subseteq Q_{R_t+|M^{-1}h|_{\max}}
$$

\noindent holds, the estimate
\begin{align*}
\|t(\tra^{-h}\xi)-t(\tra^{-h}\eta)\|_{\rm op} \;
	\leq\; &C_t \,  \max_{x\in Q_{R_t}+h\cap\Ll}\, d_\as\big(\xi(x),\eta(x)\big)^\beta\\
	\leq\; &C_t \,  \max_{x\in Q_{R_t+|M^{-1}h|_{\max}}\cap\Ll}\, d_\as\big(\xi(x),\eta(x)\big)^\beta
\end{align*}

\noindent follows. Hence, assertion (b) is proven.
\hfill $\Box$

\vspace{.2cm}

\noindent Consider the operator family $H=(H_\xi)_{\xi\in\csp}\in\CG$ of the form given in Equation~\eqref{eq:Ham} where $t_h:\csp\to M_N(\CM)$ is strongly pattern equivariant with $\beta$-H\"older continuous coefficients for each $h\in\Rr$. Then assertion (a) and (b) of Lemma~\ref{LipZd.Lem-BasHam} imply that $H$ is self-adjoint if

\begin{itemize}
\item[\textbf{(R1)}] $\Rr=-\Rr$;
\item[\textbf{(R2)}] the function $t_h$ satisfies $t_{-h}(\xi)=t_h^\ast(\tra^{-h}\xi)$ where $t_h^\ast(\tra^{-h}\xi)$ denotes the adjoint of the matrix $t_h(\tra^{-h}\xi)$.
\end{itemize}

\noindent Lemma~\ref{LipZd.Lem-BasHam}~(b) implies that the function $t_{-h}$ in \textbf{(R2)} is still strongly pattern equivariant with $\beta$-H\"older continuous coefficients with the same constant $C_{t_h}$ but different radius of influence. 

\vspace{.1cm}

\noindent For $\beta\geq 0$ , the {\em Schur $\beta$-norm} of $H$ was already defined in \eqref{eq:Schur} and denoted by $\|H\|_\beta$. It is straightforward to check that $\|H\|\leq\|H\|_\beta$ holds for all $\beta\geq 0$. Thus, if $H=(H_\xi)_{\xi\in\csp}$ is given by Equation~\eqref{eq:Ham} with infinite range $\Rr$ and $\|H\|_\beta<\infty$ for some $\beta\geq 0$, then $H\in\CG$ follows. In this case $H$ is approximated in the $C^\ast$-norm by the {\em restriction} $H|_r:=\big(H_\xi|_r\big)_{\xi\in\csp}$ defined by
\begin{align}\label{eq:RestHam}
H_\xi|_r \; 
	:= \; \sum_{h\in\Rr\cap Q_r} 
		\htt_{h,\xi} \, U_h \,.
\end{align}

\noindent More precisely, $\|H-H|_r\|\to 0$ follows if $r\to\infty$ and $\|H\|_\beta<\infty$. It is worth mentioning that $-Q_r=Q_r$ holds implying $-(\Rr\cap Q_r)=\Rr\cap Q_r$. Thus, the range set $\Rr_r:=\Rr\cap Q_r$ of $H|_r$ satisfies \textbf{(R1)} if $\Rr$ does so. 

\vspace{.1cm}

\noindent The spectrum of a Hamiltonian is studied in this work. Let $H=(H_\xi)_{\xi\in\csp}\in\CG$ be self-adjoint. Each operator $H_\xi:\Hh\to\Hh$ has spectrum $\sigma(H_\xi)\subseteq\RM$ being compact and non-empty. For $\Xi\in\Inv$, define the {\em spectrum of the operator family $H_\Xi:=(H_\xi)_{\xi\in\Xi}$} by
$$\sigma(H_\Xi) 
	\;:=\; \overline{\bigcup_{\xi\in\Xi} \sigma(H_\xi)} \subseteq\RM\,.
$$ 

\noindent Since $\sup_{\xi\in\Xi}\|H_\xi\|\leq\|H\|$ is finite, $\sigma(H_\Xi)$ is a non-empty compact subset of $\RM$. Furthermore, $\rho(H_\Xi):=\CM\setminus\sigma(H_\Xi)$ is called the {\em resolvent set of the operator family $H_\Xi:=(H_\xi)_{\xi\in\Xi}$}, which is an open subset of $\CM$. If $\mathfrak{h}$ is the element in the \Cs induced by the transformation groupoid $\Xi\rtimes_{\tra}\Ll$ that defines $H_\Xi$, then $\sigma(\mathfrak{h})=\sigma(H_\Xi)$, see e.g. \cite{LS03,Bec16,BBdN17}. It is worth mentioning that the closure in the definition of $\sigma(H_\Xi)$ is not necessary since the union of the spectra is already closed by the amenability of the group $\Ll\subseteq\RM^d$ \cite{Ex14,NP15}. Then the distance of two spectra $\sigma(H_\Xi)$ and $\sigma(H_\Theta)$ for $\Xi,\Theta\in\Inv$ is measured by the Hausdorff distance on the compact subset of $\RM$, namely
$$d_H\big(\sigma(H_\Xi),\sigma(H_\Theta)\big) \; := \;
	\max\left\{
		\sup_{E_1\in\sigma(H_\Xi)} 
			\inf_{E_2\in\sigma(H_\Theta)} |E_1-E_2|
		\,,\;
		\sup_{E_2\in\sigma(H_\Theta)} 
			\inf_{E_1\in\sigma(H_\Xi)} |E_1-E_2|
	\right\}
	\,.
$$ 

\noindent The subsection is finished with two auxiliary statements that are used in Proposition~\ref{LipZd.Prop-EsRes} which is the heart of the proof of the main theorem. Recall that $U_h:\Hh\to\Hh\,,\, U_h\psi=\psi(\cdot -h)\,,$ is the translation operator acting on the Hilbert space $\Hh:=\ell^2(\Ll)\otimes\CM^N$. For the sake of simplicity, the symbol $U_h$ is also used for the translation operator acting on the Hilbert space $\ell^2(\Ll)$.

\vspace{.1cm}

\noindent For $\varphi\in\ell^2(\Ll)$, we write $\varphi\geq 0$ if $\varphi(z)\geq 0$ for all $z\in\Ll$. Then a linear bounded operator $A:\ell^2(\Ll)\to\ell^2(\Ll)$ is called {\em positivity preserving}, if $\varphi\geq 0$ implies $A\varphi\geq 0$.

\begin{lemma}
\label{LipZd.Lem-Hbeta}
Let  $0< \beta\leq 1$ and $H:=\big(H_\xi:\Hh\to\Hh\big)_{\xi\in\csp}$ be a strongly pattern equivariant Hamiltonian with $\beta$-H\"older continuous coefficients. Let $t_{h,\beta}:\csp\to\CM$ be defined as $t_{h,\beta}(\xi):=(1+|h|^2)^{\beta/2}\, \|t_h(\xi)\|_{\rm op}$ and $\htt_{h,\beta,\xi}:\Ll\to\CM\,, \, x\mapsto t_{h,\beta}(\tra^{-x}\xi)$. Define $H_\xi^{\beta}:\ell^2(\Ll)\to\ell^2(\Ll)$ by:
$$H_\xi^\beta \; 
	:= \; \sum_{h\in\Rr} 
		\htt_{h,\beta,\xi}\, U_h.
$$

Then $H^\beta:=\big(H_\xi^{\beta}\big)_{\xi\in\csp}$ is a positivity preserving self-adjoint operator family satisfying \textbf{(R1)}, \textbf{(R2)} and its operator norm $\|H^\beta\|:= \sup_{\xi\in\csp} \|H_\xi^{\beta}\|$ is bounded by $\|H\|_{\beta}$.
\end{lemma}

\noindent {\bf Proof: } By definition, $H^\beta_\xi:\ell^2(\Ll)\to\ell^2(\Ll)$ is a linear operator. Then the estimate (see \eqref{eq:Schur} for the definition of $\|t\|_\infty$) 
$$\|H_\xi^\beta\| \; 
	\leq \; \sum_{h\in\Rr} 
		\left\| \htt_{h,\beta,\xi} \right\| \, \|U_h\|
	\leq \;  \sum_{h\in\Rr} \|t_h\|_\infty \, (1+h^2)^{\frac{\beta}{2}}=\|H\|_\beta
$$

follows. The range $\Rr$ of $H$ (and so of $H^\beta$) satisfies \textbf{(R1)}. Furthermore,  
$$t_{-h,\beta}(\xi)
	= (1+|h|^2)^{\beta/2}\|t_{-h}(\xi)\|_{\rm op}
	= (1+|h|^2)^{\beta/2}\|t_h^\ast(\tra^{-h}\xi)\|_{\rm op}
	= t_{h,\beta}^\ast(\tra^{-h}\xi)
$$

\noindent is derived as $t_h\,,\, h\in\Rr\,,$ satisfy \textbf{(R2)}. Thus, the functions $t_{h,\beta}:\csp\to[0,\infty)\,,\, h\in\Rr\,,$ also satisfy \textbf{(R2)} implying that $H^\beta$ is self-adjoint by Lemma~\ref{LipZd.Lem-BasHam}. Clearly, the translation operator $U_h$ is positivity preserving. Furthermore, the composition and sum of two positivity preserving operators is positivity preserving. Thus, the operator $H_\xi^\beta$ is positivity preserving for each $\xi\in\csp$ since $t_{h,\beta}\geq 0$.
\hfill $\Box$

\vspace{.2cm}

\begin{lemma}
\label{LipZd.Lem-Hinfty}
Let  $0< \beta\leq 1$ and $H:=\big(H_\xi:\Hh\to\Hh\big)_{\xi\in\csp}$ be a strongly pattern equivariant Hamiltonian with $\beta$-H\"older continuous coefficients. Define $H^\infty:\ell^2(\Ll)\to\ell^2(\Ll)$ by
$$H^\infty \; 
	:= \; \sum_{h\in\Rr} 
		\|t_h\|_\infty\, U_h\,.
$$

\noindent Then $H^\infty$ is a positivity preserving, self-adjoint operator such that $\|H^\infty\|\leq\|H\|_{\beta}$.
\end{lemma}

\noindent {\bf Proof: } The statement follows similarly as the previous one.
\hfill $\Box$


\section{Technical tools}
\label{LipZd.Chap-tech}

\noindent In this section some technical tools are introduced and auxiliary statements are proven. The first part deals with Lipschitz-partitions of unity that are used to localize the Hamiltonians as described in Section~\ref{LipZd.ssect-cometh}. Secondly, the so called superoperators are introduced, which are Lipschitz continuous maps on the Hilbert space $\Hh$ onto $\ell^2(\Ll)$. It is important to notice that these superoperators are not linear. However, for the purpose of this paper, it suffices that they are bounded.

 \subsection{Partition of unity}
 \label{LipZd.ssect-qupa}

A function $\qp:\RM^d\to\RM$ is called {\em Lipschitz continuous} if there is a constant $C>0$ such that $|\qp(x)-\qp(y)|\leq C |x-y|$ for all $x,y\in\RM^d$. The smallest constant satisfying the previous estimate is called {\em Lipschitz constant}, which is denoted by $C_L$. A family of functions $(\qp_i)_{i\in I}$ is called {\em uniformly Lipschitz continuous} if $C_L:=\sup_i C_L(i)<\infty$ where $C_L(i)>0$ is the Lipschitz constant of $\qp_i$. The set of Lipschitz continuous functions is denoted by $\Lip$. 

\vspace{.1cm}

\noindent The notion of Lipschitz-partition of unity will play a crucial role in this work. Such a partition can be chosen to be subordinate to any given covering  $(V_i)_{i\in I}$ of $\RM^d$, in general. Throughout this work the index set $I$ will be the lattice $\Ll$ and a specific covering is chosen. This cover is assumed to be uniformly locally finite, which reflects in condition \textbf{(P2)} below. Specifically, the covering $(V_i)_{i\in I}$ is called {\em uniformly locally finite} if there is an $\ns\in\NM$ such that for each $x\in\RM^d$ there are at most $\ns$ sets $V_i$ that contain $x$.

\begin{defini}
\label{LipZd.Def-LiPa}
A family of functions $(\qp_z)_{z\in\Ll}\subseteq\Lip$ with $0\leq \qp_z\leq 1$ is called a {\em Lipschitz-partition of unity} if the family $(\qp_z)_{z\in\Ll}$ is uniformly Lipschitz continuous and
\begin{itemize}
\item[\textbf{(P1)}] $\sum_{z\in\Ll} \qp_z(x) = 1$ for all $x\in\RM^d$;
\item[\textbf{(P2)}] The set $V_z:=\big\{z'\in\Ll\;:\; \supp\big(\qp_{z'}\big)\cap\supp\big(\qp_z\big)\neq\emptyset\big\}$ is finite uniformly in $z\in\Ll$, namely, $\ns:=\ns(\qp):=\sup_{z\in\Ll}\sharp V_z<\infty$.
\end{itemize}
\end{defini}

\noindent For the sake of clarity, we will construct a concrete example of such a partition.

\begin{exam}
\label{LipZd.Exam-Part} {\em For $r>0$, recall the notion of the deformed cube 
$Q_r
	:=\big\{M x\in\RM^d \,:\, |x|_{\max} \leq r \big\}.
$ 
Let $0\leq \psi\leq 1$ be Lipschitz continuous,
satisfying $\psi(x)=1$ for $x\in K:=\big\{x\in\RM^d\,:\, |x|_{\max}\leq \frac{1}{2}\big\}$ and $\supp(\psi)\subseteq U:=\big\{x\in\RM^d \,:\, |x|_{\max}< \frac{2}{3} \big\}\subseteq\RM^d$.  Since $(x+U)_{x\in\ZM^d}$ is a uniformly locally finite open cover of $\RM^d$, $\big(z+ MU\big)_{z\in\Ll}$ is so as well as $M$ is invertible where $MU:=\{Mx\,:\, x\in U\}$. Since the cover is uniformly locally finite, there is an $\ns\in\NM$ satisfying $\sharp \{z\in\Ll\,:\, x\in z+MU\}\leq \ns$ uniformly in $x\in\RM^d$. The constraint $\psi(x)=1$ for $x\in K$ implies $\sum_{z'\in\Ll} \psi\big(M^{-1}(x-z')\big)\geq 1$.
 
\vspace{.1cm}

\noindent Define the map $\qp:\RM^d\to[0,1]$ by $\qp(x):=\frac{\psi(M^{-1}x)}{\sum_{z'\in\Ll} \psi(M^{-1}(x-z'))}$. Its support $\supp(\qp)$ is contained in $MU\subseteq Q_{\frac{2}{3}}$. It is also straightforward to check that the family $(\qp_z)_{z\in\Ll}$ defined by $\qp_z(x):=\qp(x-z)$ satisfies all the conditions in Definition~\ref{LipZd.Def-LiPa}.
}
\hfill $\Box$
\end{exam}

\begin{rem}
\label{LipZd.Rem-ns3}
{\em It is worth mentioning that $\ns\geq 3$ holds for a Lipschitz-partition of unity in $\RM$. Thus, $\ns\geq 3$ follows for any Lipschitz-partition of unity in $\RM^d$. For indeed, if $\ns\leq 2$ one could induce a Lipschitz-partition of unity for $\RM$ by restricting the partition to the first component in $\RM^d$, a contradiction.
}
\hfill $\Box$
\end{rem}

\begin{defini}
\label{LipZd.Def-SpePart}
Let $r>0$ and $\qp\in\Lip$ be such that $\qp_z:=\qp(\cdot -z)\,,\, z\in\Ll\,,$ defines a Lipschitz-partition of unity with $\ns:=\ns(\qp)$, Lipschitz constant $C_L>0$ and $\supp(\qp)\subseteq Q_{1}$. Define the family of functions $\big(\qp_z^{(r)}\big)_{z\in\Ll}$ by 
$$\qp_z^{(r)}:\RM^d\to\RM\,,\qquad\qp_z^{(r)}(x):=\qp\left(\frac{x}{r}-z\right)\,.
$$
Furthermore, $\chi^{(r)}_z:\RM^d\to\{0,1\}$ denotes the characteristic function of the support $\supp\big(\qp_z^{(r)}\big)\subseteq\RM^d$.
\end{defini}

Example~\ref{LipZd.Exam-Part} shows that $\qp\in\Lip$ exists such that $(\qp_z)_{z\in\Ll}$ is a Lipschitz-partition of unity. The latter defined family $\big(\qp_z^{(r)}\big)_{z\in\Ll}$  of functions turns out to be also a Lipschitz-partition of unity with same bound in \textbf{(P2)}:

\begin{lemma}
\label{LipZd.Lem-LiPa}
Let $r>0$ and $\qp\in\Lip$ with Lipschitz constant $C_L>0$ be such that $(\qp_z)_{z\in\Ll}$ is a Lipschitz-partition of unity. Then the family $\big(\qp_z^{(r)}\big)_{z\in\Ll}$ defined in Definition~\ref{LipZd.Def-SpePart} is a Lipschitz-partition of unity with Lipschitz constant $\frac{C_L}{r}$ satisfying $\ns(\qp)=\ns\big(\qp^{(r)}\big)$. If $\supp(\qp)\subseteq Q_s$ for some $s>0$, then the support of $\qp^{(r)}_z$ is contained in $rz + Q_{sr}$.
\end{lemma}

\noindent {\bf Proof: } It is immediate to see the estimate
\begin{align}\label{eq:qp-r}
\big|\qp_z^{(r)}(x)-\qp_z^{(r)}(y)\big| 
	\leq \frac{C_L}{r} \, \big| x-y \big|\, , \qquad \; x,y\in \mathbb{R}^d.
\end{align} 

Condition \textbf{(P1)} follows by a short computation while \textbf{(P2)} and $\ns(\qp)=\ns\big(\qp^{(r)}\big)$ are derived from the identity
$$V_z \; 
	= \; \left\{\left.
		z'\in\Ll
		\;\right|\;
		\supp\big(\qp_{z'}^{(r)}\big)\cap\supp\big(\qp_z^{(r)}\big)\neq\emptyset
	\right\}\,.
$$

Finally, it is straightforward to show $\supp\big(\qp^{(r)}_z\big)\subseteq rz + Q_{sr}$.
\hfill $\Box$

\vspace{.2cm}

\begin{lemma}
\label{LipZd.Lem-LiQuPa}
Let $r>0$ and $(\qp^{(r)}_z)_{z\in\Ll}$ be the Lipschitz-partition of unity defined in Definition~\ref{LipZd.Def-SpePart}. For $0\leq\beta\leq 1$ and $z\in\Ll$, the estimate
$$\big| \qp^{(r)}_z(x) - \qp^{(r)}_z(y)\big)\big| \;\leq\;
	\left(\frac{|x-y|}{r}\right)^\beta \; 
		2^{1-\beta} \, C_L^\beta \;
			\Big( \chi^{(r)}_z(x) + \chi^{(r)}_z(y) \Big)
$$

\noindent holds for all $x,y\in\RM^d$. 
\end{lemma}

\noindent {\bf Proof: } Combining \eqref{eq:qp-r} with the inequality $\big| \qp^{(r)}_z(x) - \qp^{(r)}_z(y)\big|\leq 2$ we obtain by interpolation:
$$\big| \qp^{(r)}_z(x) - \qp^{(r)}_z(y)\big| \; 
	\leq \; \left(
		C_L \, \frac{|y-x|}{r} 
	\right)^\beta 
		2 ^{1-\beta},\quad 0\leq\beta\leq 1.
$$

The last ingredient is the identity
$$\big| \qp^{(r)}_z(x) - \qp^{(r)}_z(y)\big)\big| \;
	= \; \big| \qp^{(r)}_z(x) - \qp^{(r)}_z(y)\big| 
		\, \big( \chi^{(r)}_z(x) + \chi^{(r)}_z(y) \big),
$$

\noindent finishing the proof.
\hfill $\Box$

 \subsection{Superoperators}
 \label{LipZd.ssect-supop}

\noindent Recall that $\Hh$ denotes the Hilbert space $\ell^2(\Ll)\otimes\CM^N$ on which the Hamiltonians act. Furthermore, $\chi^{(r)}_z:\RM^d\to\{0,1\}$ denotes the characteristic function of the support $\supp(\qp^{(r)}_z)$. The corresponding multiplication operator by the function $\chi^{(r)}_z$ on $\Hh$ (acting as the identity on $\mathbb{C}^N$) is denoted by the symbol $\Mch^{(r)}_z$, with $\|\Mch_z^{(r)}\|=1$  for each $z\in\Ll$. In the following $\Cc_c(\Ll)$ denotes the set of functions $\varphi:\Ll\to\CM$ with finite support in $\Ll$. Note in the following that $(B\varphi)(x)\in\CM^N$ holds if $B:\Hh\to\Hh$, $\varphi\in\Hh$ and $x\in\Ll$. In this case $|(B\varphi)(x)|$ denotes the Euclidean length of the vector $B\varphi(x)$.

\begin{lemma}
\label{LipZd.Lem-G(A)}
Consider an operator family $(A_z)_{z\in\Ll}\subseteq\Bb\big(\Hh\big)$ with $\|A\| := \sup_{z\in\Ll} \|A_z\|< \infty$. Then the map $\OG(A):\Cc_c(\Ll)\otimes\CM^N\to\ell^2(\Ll)$ defined by
$$\big(\OG(A)\varphi\big) (x)\; 
	:= \; \sum_{z\in\Ll} \Mch^{(r)}_z(x)\, \big|\big (A_z \, \Mch^{(r)}_z \varphi\big)(x)\big|\,,
	\qquad
	\varphi\in\Cc_c(\Ll)\otimes\CM^N\,,\, x\in \Ll,
$$

\noindent satisfies 
$$\big\|\OG(A)\varphi_1 - \OG(A)\varphi_2 \big\|_\Hh
	\leq \ns \, \|A\| \, \|\varphi_1-\varphi_2\|_\Hh\,,
	\qquad \varphi_1,\varphi_2\in\Cc_c(\Ll)\otimes\CM^N\,.
$$

\noindent Furthermore, $\OG(A)$ uniquely extends to a continuous bounded map on $\Hh$ to $\ell^2(\Ll)$ such that 
$$\sup\big\{ \|\OG(A)\varphi\|_\Hh \;:\; \varphi\in\Hh \, \text{ with } \, \|\varphi\|_\Hh\leq 1\big\}
		\;\leq\; \ns \, \|A\|\,.
$$
\end{lemma}

\noindent {\bf Proof: } As introduced in Definition~\ref{LipZd.Def-LiPa}, $V^{(r)}_z$ denotes the set of all $z'\in\Ll$ such that $\supp(\qp^{(r)}_z)\cap\supp(\qp^{(r)}_z)\neq\emptyset$. Let $\varphi\in\Cc_c(\Ll)\otimes\CM^N$. First note that 
$$
\varphi_z(x):=\big|\big (A_z \, \Mch^{(r)}_z \varphi\big)(x)\big|
	= \left(
		\sum_{k=1}^N \big|\big (A_z \, \Mch^{(r)}_z \varphi\big)_k(x)\big|^2
	\right)^{\frac{1}{2}}
	\,,\qquad x,z\in\Ll\,,
$$

\noindent defines an element in $\ell^2(\Ll)$. Thus, the Cauchy-Schwarz inequality on $\ell^2(\Ll)$ yields
$$
\big|\big\langle
	\Mch^{(r)}_z\, \varphi_z\,,\; \Mch^{(r)}_{z'}\, \varphi_{z'}
\big\rangle_{\ell^2(\Ll)}\big|
	\leq \big\|\Mch^{(r)}_z\, \varphi_z\big\|_{\ell^2(\Ll)}\; \big\|\Mch^{(r)}_{z'}\, \varphi_{z'}\big\|_{\ell^2(\Ll)}
	\leq \|\varphi_z\|_{\ell^2(\Ll)}\; \|\varphi_{z'}\|_{\ell^2(\Ll)}\,.
$$

\noindent Note that the latter inner product vanishes if $z\in\Ll$ and $z'\not\in V^{(r)}_z$. Furthermore, a short computation leads to
$$
\|\varphi_z\|_{\ell^2(\Ll)}^2 
	= \big\| A_z \Mch^{(r)}_z\varphi\big\|_\Hh^2 
	\leq \|A\|^2 \| \Mch^{(r)}_z\varphi\|_\Hh^2\,.
$$

\noindent Since $2ab\leq a^2+b^2$ for $a,b\geq 0$, the previous considerations yield
\begin{align*}
\|\OG(A)\varphi\|^2 \;
	\leq \; &\sum_{z\in\Ll}\, \sum_{z'\in V^{(r)}_{z}} 
		\big|
			\big\langle
				\Mch^{(r)}_z\, |\big (A_z \, \Mch^{(r)}_z \varphi\big )|\,,\;
				\Mch^{(r)}_{z'}\, |\big (A_{z'} \, \Mch^{(r)}_{z'} \varphi\big )|
			\big\rangle_{\ell^2(\Ll)}
		\big|\\
	\leq \; &\frac{\|A\|^2}{2} \, 
		\sum_{z\in\Ll}\, \sum_{z'\in V^{(r)}_{z}} 
			\big(
				\|\Mch^{(r)}_z \varphi\|_\Hh^2 + \|\Mch^{(r)}_{z'} \varphi\|_\Hh^2
			\big)\,.
\end{align*}

\noindent Let $l_x$ be the number of $z'\in\Ll$ satisfying that $\chi^{(r)}_{z'}(x)\neq 0$ for fixed $x\in\Ll$. Thus, $l_x\leq \ns$ is derived for all $x\in\Ll$ by \textbf{(P2)} and Lemma~\ref{LipZd.Lem-LiPa}. Hence, 
$$\sum_{z\in\Ll}\, \sum_{z'\in V^{(r)}_{z}} \chi^{(r)}_{z'}(x) \;
	\leq\; \ns \, \sharp\big\{z'\in V^{(r)}_{z} \,:\, \chi^{(r)}_{z'}(x)\neq 0\big\}
	\leq\; \ns^2
$$

\noindent follows. Using the previous considerations, the estimate
$$\sum_{z\in\Ll}\, \sum_{z'\in V^{(r)}_{z}} \|\Mch^{(r)}_{z'}\varphi\|_\Hh^2 \;
	= \; \sum_{x\in\Ll}\sum_{k=1}^N |\varphi_k(x)|^2 \
		\left(
			\sum_{z\in\Ll}\, \sum_{z'\in V^{(r)}_{z}} \chi^{(r)}_{z'}(x)
		\right)\\
	\leq \; \|\varphi\|_\Hh^2 \, \ns^2
$$

\noindent is proven. Similarly,
$$\sum_{z\in\Ll}\, \sum_{z'\in V^{(r)}_{z}} \|\Mch^{(r)}_z\varphi\|_\Hh^2 \;
	= \; \sum_{x\in\Ll}\sum_{k=1}^N |\varphi_k(x)|^2 \, \ns \, 
		\sharp\left\{
			z\in\Ll\,\big|\,
			\chi^{(r)}_z(x)\neq 0
		\right\}
	\leq \; \|\varphi\|_\Hh^2 \, \ns^2
$$

\noindent is deduced. These observations imply $\|\OG(A)\varphi\|_\Hh^2 \leq \ns^2 \, \|A\|^2 \, \|\varphi\|_\Hh^2$ for all $\varphi\in\Cc_c(\Ll)$. Let $\varphi_1,\varphi_2\in\Cc_c(\Ll)\otimes\CM^N$. Since we have:
\begin{align*}
\big| \big (\OG(A)\varphi_1-\OG(A)\varphi_2\big )_k (x)\big|
	&\leq \sum_{z\in\Ll} 
		\Mch^{(r)}_z \, 
		\big|
			|\big (A_z\, \Mch^{(r)}_z \varphi_1\big )_k(x)| - |\big (A_z\, \Mch^{(r)}_z \varphi_2\big )_k(x)|
		\big| \\
		&	\leq \sum_{z\in\Ll} 
		\Mch^{(r)}_z \,
		\big|
			\big (A_z\, \Mch^{(r)}_z (\varphi_1-\varphi_2)\big )_k(x)
		\big|\,,
\end{align*}

the previous considerations lead to
$$\big\| \OG(A)\varphi_1-\OG(A)\varphi_2\big\|_\Hh
	\leq \big\|\OG(A)(\varphi_1-\varphi_2)\big\|_\Hh
	\leq \ns \, \|A\| \, \|\varphi_1-\varphi_2\|_\Hh\,,
$$

\noindent namely $\OG(A)$ is Lipschitz continuous. Thus, there is a unique continuous extension $\OG(A):\Hh\to\ell^2(\Ll)$ satisfying
$$\sup\big\{ \|\OG(A)\varphi\|_\Hh \;:\; \varphi\in\Hh \, \text{ with } \, \|\varphi\|_\Hh\leq 1\big\}
		\;\leq\; \ns \, \|A\|
$$

\noindent as $\Cc_c(\Ll)\otimes\CM^N\subseteq\Hh$ is dense.
\hfill $\Box$

\vspace{.2cm}

\begin{lemma}
\label{LipZd.Lem-G_D(A)}
Let $B:\ell^2(\Ll)\to\ell^2(\Ll)$ be a positivity preserving, linear bounded operator. Consider an operator family $(A_z)_{z\in\Ll}\subseteq\Bb(\Hh)$ such that $\|A\| := \sup_{z\in\Ll} \|A_z\| < \infty$. Then the map $\OG_B(A):\Cc_c(\Ll)\otimes\CM^N\to\ell^2(\Ll)$ defined by
$$\big(\OG_B(A)\varphi\big)(x) \; 
	:= \; \sum_{z\in\Ll} \Mch^{(r)}_z(x)\, \left (B\, \big|A_z\, \Mch^{(r)}_z \varphi\big|\right )(x)\,,
	\qquad
	\varphi\in\Cc_c(\Ll)\otimes\CM^N\,,\, x\in \Ll,
$$

\noindent satisfies 
$$\big\|\OG_B(A)\varphi_1 - \OG_B(A)\varphi_2 \big\|_\Hh
	\leq \ns \, \|A\| \, \|B\| \, \|\varphi_1-\varphi_2\|_\Hh\,,
	\qquad \varphi_1,\varphi_2\in\Cc_c(\Ll)\otimes\CM^N\,.
$$

\noindent Furthermore, $\OG_B(A)$ extends to a continuous bounded map on $\Hh$ to $\ell^2(\Ll)$ such that
$$\sup\big\{ \|\OG_B(A)\varphi\|_\Hh \;:\; \varphi\in\Hh \, \text{ with } \, \|\varphi\|_\Hh\leq 1\big\} \;
	\leq \; \ns \, \|A\| \, \|B\|\,.
$$
\end{lemma}

\noindent {\bf Proof: } Recall that $\big|A_z\, \Mch^{(r)}_z \varphi\big|$ is an element of $\ell^2(\Ll)$ and so $B\big|A_z\, \Mch^{(r)}_z \varphi\big|$ is well-defined. Following the lines of the proof Lemma~\ref{LipZd.Lem-G(A)} we obtain
$$\|\OG_B(A)\varphi\|_\Hh \;
	\leq \; \ns \, \|B\| \, \|A\| \, \|\varphi\|_\Hh\,.
$$

Since $B$ is a positivity preserving linear operator we have
$$0 \leq  B\big( |\varphi-\psi| - |\varphi|+|\psi|\big) 
	= B|\varphi-\psi| - \big( B|\varphi| - B|\psi|\big)\,.
$$

This implies a pointwise estimate:
$$\Big|B\big|A_z\, \Mch^{(r)}_z \varphi_1\big| - B\big|A_z\, \Mch^{(r)}_z \varphi_2\big|\Big| 
	\leq B\big|A_z\, \Mch^{(r)}_z (\varphi_1-\varphi_2)\big|\,.
$$

The latter yields the pointwise estimate: 
$$
\big| \OG_B(A)\varphi_1-\OG_B(A)\varphi_2 \big| \; 
	\leq \; \sum_{z\in\Ll} 
		\Mch^{(r)}_{z} \, B\,
		\big|
			A_z\, \Mch^{(r)}_z (\varphi_1-\varphi_2)
		\big| \;
	= \; \OG_B(A)(\varphi_1-\varphi_2)
$$

which leads to
$$
\big\|\OG_B(A)\varphi_1 - \OG_B(A)\varphi_2\big\|_\Hh \;
	\leq \; \big\|\OG_B(\varphi_1 - \varphi_2)\big\|_\Hh
	\leq \; \ns \, \|B\| \, \|A\| \, \|\varphi_1-\varphi_2\|_\Hh \,.
$$

Then the density of $\Cc_c(\Ll)\otimes\CM^N\subseteq\Hh$ finishes the proof.
\hfill $\Box$


\section{Proofs of the main results}
\label{LipZd.Chap-proof}

\noindent The following two auxiliary lemmas provide estimates on the potential and the kinetic term if commuted with the Lipschitz partition of unity.

\begin{lemma}
\label{LipZd.Lem-ComPo}
For each $r>1$, let $(\qp_z^{(r)})_{z\in\Ll}$ be the Lipschitz-partition of unity defined in Definition~\ref{LipZd.Def-SpePart}. Let $0<\beta\leq 1$ and $t:\csp\to M_N(\CM)$ be a strongly pattern equivariant function with $\beta$-H\"older continuous coefficients, constant $C_t\geq 1$ and radius of influence $R_t\geq 1$. Let $\Xi,\Theta\in\Inv$ be such that 
$$
r:=\Hdco(\Xi,\Theta)^{-1} 
	> \tR:= R_t+\|M^{-1}\|_{\max} \, \|M\|_{\max} +1
$$

Then for each $\xi\in\Xi$ and $z\in\Ll$, there exists an $\eta_\xi:=\eta(\xi,z,r,\widetilde{R},\Ll)\in\Theta$ such that:
$$\Big|\Big(\htt_{\xi}\, \Mqp^{(r-\tR)} - \Mqp^{(r-\tR)}\, \htt_{\eta_z}\Big)\varphi\Big|(x)\;
	\leq \; \frac{2}{r^\beta} \, C_t \, \chi_z^{(r-\tR)}(x) \, \big|\varphi(x)\big|
	\,,\qquad
	\varphi\in\Hh\,,\; x\in\Ll\,,
$$

\noindent holds.
\end{lemma}

\noindent {\bf Proof: } Let $z\in\Ll$. For each $x\in\RM^d$ we have:
$$\min_{y\in\Ll} |x-y|_{\max}
	\;\leq\; \|M\|_{\max}  \min_{m\in\Zd} |M^{-1}x-m|_{\max} 
	\;\leq\; \|M\|_{\max} 
$$

Thus there exists an $z_r\in\Ll$ such that $\big|z_r-(r-\tR)\, z \big|_{\max}\leq\|M\|_{\max}$. Since $\Xi$ is invariant, we have $\tra^{-z_r}\xi\in \Xi$ and due to  Lemma~\ref{LipZd.Lem-BaEq} we may find $\tilde{\eta}_\xi\in\Theta$ such that 
$$\dco\big(\tra^{-z_r}\xi, \tilde{\eta}_\xi\big)
	\;\leq\; \Hdco(\Xi,\Theta)
	\,.
$$

\noindent By definition, $R_t\geq 1$ holds implying $r>2$. Since $\Theta$ is invariant, $\eta_\xi:=\tra^{z_r}\tilde{\eta}_\xi\in \Theta$.  Then Lemma~\ref{LipZd.Lem-BaEq}~(b) implies
\begin{equation}\label{eq:PotMet}
d_\as\big(\tra^{-z_r}\xi(y),\tra^{-z_r}\eta_\xi(y) \big)
	= d_\as\big(\tra^{-z_r}\xi(y), \tilde{\eta}_\xi(y)\big)
	\leq \frac{1}{r-1} \leq \frac{2}{r}
	\,, \qquad
	y\in Q_{r-1}\cap\Ll\,.
\end{equation}

\noindent Consider some $\varphi\in\Hh$ and $x\in\Ll$. Then a short computation leads to
$$\Big(\htt_{\xi}\, \Mqp^{(r-\tR)}\varphi - \Mqp^{(r-\tR)}\, \htt_{\eta_\xi}\varphi\Big)(x)
	= \qp_z^{(r-\tR)}(x)\;\Big(t\big(\tra^{-x}\xi\big)\varphi - t\big(\tra^{-x}\eta_\xi\big)\varphi \Big)\;  (x)\,.
$$

\noindent Clearly, the term vanishes if $x\not\in\supp\big(\qp^{(r-\tR)}_z\big)$ and $0\leq \qp_z^{(r-\tR)}(x)\leq \chi_z^{(r-\tR)}(x)$. Thus, it suffices to show
$$\|t\big(\tra^{-x}\xi\big) - t\big(\tra^{-x}\eta_\xi\big)\|_{\rm op}
	\;\leq\; 2\, C_t \,  r^{-\beta}
	\,,\qquad x\in\supp\big(\qp^{(r-\tR)}_z\big)\cap\Ll\,.
$$

\noindent Let $x\in\supp\big(\qp^{(r-\tR)}_z\big)\cap\Ll$. According to Definition~\ref{LipZd.Def-SpePart}, the inclusion $\supp(\qp)\subseteq Q_1$ holds. Thus, Lemma~\ref{LipZd.Lem-LiPa} asserts $x\in \supp\big(\qp^{(r-\tR)}_z\big)\subseteq(r-\tR)\, z + Q_{r-\tR}$. 

\vspace{.1cm}

Using the definition of $\tR$, we have $|x_1+x_2|_{\max}\leq r-\|M^{-1}\|_{\max}\,\|M\|_{\max}-1$ whenever $|x_1|_{\max}\leq r-\tR$ and $|x_2|_{\max}\leq R_t$. Thus, $Q_{r-\tR} + Q_{R_t}$ is contained in $Q_{r-\|M^{-1}\|_{\max} \, \|M\|_{\max}-1}$ implying 
$$x + Q_{R_t}
	\subseteq (r-\tR)\, z + Q_{r-\|M^{-1}\|_{\max} \, \|M\|_{\max}-1}
		=  z_r + \left((r-\tR)\, z - z_r\right)  + Q_{r-\|M^{-1}\|_{\max} \, \|M\|_{\max}-1}\,.
$$

\noindent In addition, $\big((r-\tR)\, z - z_r\big)  + Q_{r-\|M^{-1}\|_{\max} \|M\|_{\max}-1}\subseteq Q_{r-1}$ is an immediate consequence of the estimates
$$\left|M^{-1}\left((r-\tR)\, z - z_r\right)+y \right|_{\max}
	\;\leq\; |y|_{\max} + \|M^{-1}\|_{\max} \, \left|(r-\tR)\, z - z_r\right|_{\max}
	\;\leq\;  r-1
$$ 

\noindent for each $y\in\RM^d$ satisfying $My\in Q_{r-\|M^{-1}\|_{\max} \|M\|_{\max}-1}$. Consequently we have:
$$x + Q_{R_t} \;
	\subseteq z_r + Q_{r-1}
	\,,\qquad x\in\supp\big(\qp^{(r-\tR)}_z\big)\cap\Ll\,.
$$ 

Since $t$ satisfies \eqref{eq:PatEq}, $x\in\supp\big(\qp^{(r-\tR)}_z\big)\cap\Ll\subseteq(r-\tR)\, z + Q_{r-\tR}\cap\Ll$ yields
\begin{align*}
\|t\big(\tra^{-x}\xi\big) - t\big(\tra^{-x}\eta_\xi\big)\|_{\rm op}
	\;&\leq\; C_t \, \max\left\{
		d_\as\big(\xi(y),\eta_\xi(y) \big)^\beta 
		\;:\; 
		y\in x+Q_{R_t}\cap\Ll
	\right\}\\
	\;&\leq\; C_t \,  \max\left\{
		d_\as\big(\tra^{-z_r}\xi(y),\tra^{-z_r}\eta_\xi(y) \big)^\beta
		\;:\;
		y\in Q_{r-1}\cap\Ll
	\right\}.
\end{align*}

Then \eqref{eq:PotMet} yields the desired estimate.
\hfill $\Box$

\vspace{.2cm}

For two operators $A,B\in\Bb(\Hh)$ we denote their commutator $AB-BA$ by $[A,B]$. The following estimate is based on Lemma~\ref{LipZd.Lem-LiQuPa}.

\begin{lemma}
\label{LipZd.Lem-ComTrans}
Let $(\qp_z^{(r)})_{z\in\Ll}$ be the Lipschitz-partition of unity defined in Definition~\ref{LipZd.Def-SpePart}. For each $z\,,\; h\in\Ll\,,\; 0\leq \beta\leq 1$ and $r\geq 1$, the estimate
$$\left|\left(\Big[U_h\,,\Mqp^{(r)} \Big] \varphi\right)(x)\right|
	\leq \frac{|h|^\beta}{r^\beta} \, 
			2^{1-\beta} \, C_L^\beta \, 
				\Big( \chi^{(r)}_z(x-h) + \chi^{(r)}_z(x) \Big)\
					\left|\big(U_h\varphi\big)(x)\right|
$$

\noindent holds for all $\varphi\in\Hh$ and $x\in\Ll$.
\end{lemma}

\noindent {\bf Proof: } Let $\varphi\in\Hh$ and $x\in\Ll$. A short computation leads to
$$\left(\Big[U_h\,,\Mqp^{(r)} \Big] \varphi\right)(x) \;
	= \;
		\left(
			\qp^{(r)}_z(x-h)-\qp^{(r)}_z(x)		
		\right)
		\; \big(U_h\varphi\big)(x)\,.
$$

\noindent In addition, Lemma~\ref{LipZd.Lem-LiQuPa} implies
$$\left| \qp^{(r)}_z(x-h)-\qp^{(r)}_z(x) \right|
	\leq \frac{|h|^\beta}{r^\beta} \, 
			2^{1-\beta} \, C_L^\beta \, 
				\Big( \chi^{(r)}_z(x-h) + \chi^{(r)}_z(x) \Big)\,.
$$
\hfill $\Box$

\vspace{.2cm}

We are now interested in constructing an approximate inverse of $H_\xi-z$ using the resolvent of $H_\eta$. Denote the distance of $\zz\in\CM$ to a compact subset $K\subseteq\CM$ by
$$\dist(\zz,K):=
	\inf\big\{
		|\zz-x|\;:\; x\in K
	\big\}\,.
$$

\noindent Recall the notion of the spectrum $\sigma(H_\Theta)$ and resolvent set $\rho(H_\Theta)$ for an $H\in\CG$ and a subshift $\Theta\in\Inv$, which were defined in Section~\ref{LipZd.ssect-ham}.

\begin{lemma}
\label{LipZd.Lem-PseuInvBoun}
Let $r>0$, $\Theta\in\Inv$ and $(\qp^{(r)}_z)_{z\in\Ll}$ be the Lipschitz-partition of unity defined in Definition~\ref{LipZd.Def-SpePart} with $\ns:=\ns\big(\qp^{(r)}\big)$ independent of $r$. Suppose $H=(H_\xi)_{\xi\in\csp}\in\CG$ is self-adjoint.  For every $z\in\Ll$, choose an arbitrary $\eta_z\in\Theta$. Then for every  $\zz\in\rho(H_\Theta)=\CM\setminus\sigma(H_\Theta)$, the operator $S(\zz)\in\Bb(\Hh)$ given by
$$S(\zz) \; 
	:= \; \sum_{z\in\Ll}
		\Mqp^{(r)} \; (H_{\eta_z}-\zz)^{-1} \; \Mch_z^{(r)}\,.
$$

\noindent is well-defined and its operator norm is bounded by $\frac{\ns}{\dist(\zz,\sigma(H_\Theta))}$.
\end{lemma}

\noindent {\bf Proof: } Note that $H_\eta-\zz$ is invertible for each $\eta\in\Theta$ as $\zz\in\rho(H_\Theta)$. Consider the operator family $A_z:=\Mqp^{(r)}\, (H_{\eta_z}-\zz)^{-1}$ for $z\in\Ll$. Its operator norm is bounded by 
\begin{align*}
\|A\| \;
	= \; \sup_{z\in\Ll} \|A_z\| \;
			\leq \; \sup_{\eta\in\Theta}\|(H_\eta-\zz)^{-1}\| \;
				= \; \sup_{\eta\in\Theta}\; \frac{1}{\dist\big(\zz,\sigma(H_\eta)\big)} \;
					= \; \frac{1}{\dist\big(\zz,\sigma(H_\Theta)\big)}
\end{align*}

Let $\varphi\in\Cc_c(\Ll)\otimes\CM^N$ and $x\in\Ll$. Since $\Mqp^{(r)}=\Mch_z^{(r)}\, \Mqp^{(r)}$, the estimate
$$\big|\big(S(\zz)\varphi\big)(x)\big| \; 
	= \; \left| \sum_{z\in\Ll}
		 \left(\Mch^{(r)}_z\, A_z \, \Mch^{(r)}_z\varphi\right)(x)\;
		\right|
	\leq\; \big(\OG(A)\varphi\big)(x)
$$

\noindent follows where $\OG(A)$ is the map defined in Lemma~\ref{LipZd.Lem-G(A)}. Hence, Lemma~\ref{LipZd.Lem-G(A)} implies $\|\big(S(\zz)\big)\| \leq \ns \, \|A\|$ which coupled with the estimate on $\|A\|$ ends the proof.
\hfill $\Box$

\vspace{.2cm}


\begin{proposi}
\label{LipZd.Prop-EsRes}
Consider a strongly pattern equivariant Hamiltonian $H:=\big(H_\xi\big)_{\xi\in\csp}$ with $\beta$-H\"older continuous coefficients of finite range with radius of influence $R_H$. Let $\Xi,\Theta\in\Inv$ be such that 
$$
r:=\Hdco(\Xi,\Theta)^{-1} 
	> \tR_H:=R_H+\|M^{-1}\|_{\max} \, \|M\|_{\max}+1
$$

Let $(\qp^{(r)}_z)_{z\in\Ll}$ be a Lipschitz-partition of unity defined in Definition~\ref{LipZd.Def-SpePart} with $\ns:=\ns\big(\qp^{(r)}\big)$ independent of $r$.  Then if $\zz\in\rho(H_\Theta)$ satisfies 
\begin{align}\label{eq:EsRes}
\dist\big(\zz,\sigma(H_\Theta)\big) \;
	> \; 
		\frac{16\, \ns \, \max\big\{ C_L \,,\, 1\big\}}{2^\beta\, (r-\tR_H)^\beta} \
		 \, C_{hop} \, \|H\|_\beta
	\,,
\end{align}

we have that $\zz\in\rho(H_\Xi)$.
\end{proposi}

\noindent {\bf Proof: } For $\xi'\in\csp$, the Hamiltonian is defined by
$$H_{\xi'}:=\sum_{h\in\Rr} \htt_{h,\xi'} \, U_h
$$

\noindent satisfying \textbf{(R1)} and \textbf{(R2)} where $t_h:\csp\to M_N(\CM)\,,\, h\in\Rr\,,$ are strongly pattern equivariant with $\beta$-H\"older continuous coefficients and $\Rr$ is finite. Furthermore, $C_{hop}=\sup_{h\in\Rr} C_{t_h}$ is finite and $1\leq R_{t_h}\leq R_H$ holds for all $h\in\Rr$. 

\vspace{.1cm}

\noindent Let $\zz\in\rho(H_\Theta)$ obeying \eqref{eq:EsRes}. In particular, $\zz\in\rho(H_\eta)$ for all $\eta\in\Theta$ and $\dist\big(\zz,\sigma(H_\eta)\big) \geq \dist\big(\zz,\sigma(H_\Theta)\big)$. It suffices to prove that $\zz\in\rho(H_\xi)$ uniformly in $\xi\in\Xi$, i.e. that there exists $\delta>0$ independent of $\xi\in\Xi$ such that $\dist\big(\zz,\sigma(H_\xi)\big)>\delta$. This implies $\zz\in\rho(H_\Xi)=\CM\setminus \overline{\bigcup_{\xi\in\Xi}\sigma(H_\xi)}$.

\vspace{.1cm}

\noindent In light of this, let $\xi\in\Xi$ and we will prove that $\zz\in\rho(H_\Xi)$ going through the following steps:
\begin{itemize}
\item[(i)] An operator $S(\zz)$ is constructed such that $(H_\xi-\zz)S(\zz)=\Id + T_1(\zz) + T_2(\zz)$ where the error terms $T_1(\zz)$ and $T_2(\zz)$ come from the kinetic and the potential terms, respectively.
\item[(ii)] It is shown that $\|T_1(\zz)\|\leq \frac{1}{4}$.
\item[(iii)] It is shown that $\|T_2(\zz)\|\leq \frac{1}{4}$.
\item[(iv)] Using (i)-(iii), $\zz\in\rho(H_\xi)$ is verified.
\end{itemize}

\vspace{.1cm}

\noindent (i): Recall that $R_{t_h}\geq 1$ holds by definition. Since $0<\Hdco(\Xi,\Theta)=\frac{1}{r}$ with $r>\tR_H\geq R_{t_h}+\|M^{-1}\|_{\max}\, \|M\|_{\max}+1\geq 2$ and $C_{t_h}\leq C_{hop}$ for all $h\in\Rr$, Lemma~\ref{LipZd.Lem-ComPo} applies with $\widetilde{R}$ replaced by $\tR_H$ and $C_{t}$ replaced by $C_{hop}$. This implies that given any $z\in\Ll$, there exists a $\eta_z=\eta(\xi,z,r,\tR_H,\Ll)\in\Theta$ satisfying
\begin{equation}
\label{LipZd.Eq-EsRes}
\left|\left(\left(
		\htt_{h,\xi}\, \Mqp^{(r-\tR_H)} - \Mqp^{(r-\tR_H)}\, \htt_{h,\eta_z}
	\right) \varphi\right)(x)\right|\;
	\leq \; \frac{2}{r^{\beta}} \, C_{hop}\, \chi_z^{(r-\tR_H)}(x)\, \big|\varphi(x)\big|
\end{equation}

\noindent for all $\varphi\in\Hh$ and $x\in\Ll$. It is worth noticing that $\eta_z$ is independent of $h\in\Rr$ by Lemma~\ref{LipZd.Lem-ComPo}. With this chosen $\eta_z\in\Theta$ for $z\in\Ll$, define $S(\zz)\in\Bb(\Hh)$ by
$$S(\zz) \; 
	:= \; \sum_{z\in\Ll}
		\Mqp^{(r-\tR_H)} \, (H_{\eta_z}-\zz)^{-1} \, \Mch_z^{(r-\tR_H)}\,.
$$

\noindent According to Lemma~\ref{LipZd.Lem-PseuInvBoun}, $S(\zz)$ is a well-defined operator and 
\begin{align}\label{eq:Est-S}
\|S(\zz)\| 
	\;\leq \; \frac{\ns}{\dist\big(\zz,\sigma(H_\Theta)\big)}\,.
\end{align}

\noindent In the following we investigate the operator product $H_\xi \, S(\zz)$. In order to shorten notation, define:
$$E_{h,z}\;
	:=\; \htt_{h,\xi}\, \Mqp^{(r-\tR_H)} - \Mqp^{(r-\tR_H)}\, \htt_{h,\eta_z}
	\,,\qquad h\in\Rr\,,\, z\in\Ll\,.
$$

\noindent Then for each $z\in\Ll$ we have
\begin{align*}
H_\xi \, \Mqp^{(r-\tR_H)} \;
	= \; &\left(\sum_{h\in\Rr} 
			\htt_{h,\xi} \, \Big[U_h \, , \Mqp^{(r-\tR_H)}\Big]
			+
			\sum_{h\in\Rr}
			E_{h,z} \, U_h
		\right)
	 + \Mqp^{(r-\tR_H)} \, H_{\eta_z}\,.
\end{align*}

Define the operators $T_1(\zz),T_2(\zz)\in\Bb(\Hh)$ by
$$T_1(\zz) 
	:=  \sum_{z\in\Ll} \sum_{h\in\Rr} 
			\htt_{h,\xi} \, \Big[U_h \, , \Mqp^{(r-\tR_H)}\Big] \, A_z(\zz)\, \Mch_z^{(r-\tR_H)}
	\,,\quad
T_2(\zz)
	:= \sum_{z\in\Ll} \sum_{h\in\Rr}
			E_{h,z} \, U_h \, A_z(\zz) \, \Mch_z^{(r-\tR_H)}\,,
$$

\noindent where $A_z(\zz):=(H_{\eta_z}-\zz)^{-1}\in\Bb(\Hh)$ for $z\in\Ll$. Then the operator norm of this operator family $(A_z(\zz))_{z\in\Ll}$ satisfies
$$\|A(\zz)\|\; 
	:= \; \sup_{z\in\Ll}\|A_z(\zz)\|\;
	\leq \; \|(H_{\eta_z}-\zz)^{-1}\| \;
		= \;\frac{1}{\dist\big(\zz,\sigma(H_{\eta_z})\big)} \;
		\leq \;\frac{1}{\dist\big(\zz,\sigma(H_\Theta)\big)}\,.
$$

\noindent According to Lemma~\ref{LipZd.Lem-LiPa}, $\big(\qp_z^{(r-\tR_H)}\big)_{z\in\Ll}$ is a Lipschitz-partition of unity satisfying \textbf{(P1)} and \textbf{(P2)}. Since $\Mqp^{(r-\tR_H)}=\Mqp^{(r-\tR_H)}\Mch_z^{(r-\tR_H)}$, \textbf{(P1)} implies that $\sum_{z\in\Ll} \Mqp^{(r-\tR_H)} \, \Mch_z^{(r-\tR_H)}$ is equal to the identity operator $\Id$. With this at hand, the previous considerations lead to
\begin{align*}
(H_\xi-\zz)S(\zz) \;
	= \; & \sum_{z\in\Ll}  
		\Mqp^{(r-\tR_H)}
			\big(H_{\eta_z}-\zz\big)\, (H_{\eta_z}-\zz)^{-1}\,
				\Mch_z^{(r-\tR_H)}
		\,+\, T_1(\zz) + T_2(\zz)\\
	= \; &\Id + T_1(\zz) + T_2(\zz)\,.
\end{align*}

\noindent (ii): Let $\varphi\in\Hh$ and $x\in\Ll$. Lemma~\ref{LipZd.Lem-ComTrans} implies
\begin{eqnarray*}
\big|T_1(\zz)\varphi(x)\big|
 &\leq&
 \sum_{z\in\Ll}
  \sum_{h\in\Rr}
    \big\|t_h\big(\tra^{-x}\xi\big)\big\|_{\rm op} \,
     \Big|\left(
     	\Big[U_h \, , \Mqp^{(r-\tR_H)}\Big] \,
      	A_z(\zz)\, \Mch^{(r-\tR_H)}_z \varphi
      \right)(x)\Big|\\
&\leq&
 \frac{2^{1-\beta}\, C_L^\beta}{(r-\tR_H)^\beta} \,
  \Bigg(
   \sum_{z\in\Ll}
    \chi^{(r-\tR_H)}_z(x) \
     \sum_{h\in\Rr} |h|^\beta \, 
      \big\|t_h\big(\tra^{-x}\xi\big)\big\|_{\rm op} \,
       \Big|\Big( 
        U_h \, A_z(\zz)\, \Mch^{(r-\tR_H)}_{z} \varphi
       \Big)(x) \Big|\\
&&\quad\qquad\qquad +
 \sum_{h\in\Rr} |h|^\beta \, 
  \big\|t_h\big(\tra^{-x}\xi\big)\big\|_{\rm op} \,
   \sum_{z\in\Ll} 
    \Big(
     \Mch^{(r-\tR_H)}_z \, \big| 
    	A_z(\zz)\, \Mch^{(r-\tR_H)}_z \varphi
      \big|
    \Big)(x-h) 
   \Bigg)\,.
\end{eqnarray*}

\noindent  Recall the notion of the Hamiltonian $H^\beta=(H_\xi^\beta)_{\xi\in\as^\Zd}$ introduced in Lemma~\ref{LipZd.Lem-Hbeta} where the hopping terms are given $t_{h,\beta,\xi}:\Ll\to[0,\infty)\,,\,x\mapsto (1+|h|^2)^{\beta/2} \, \|t_h(\tra^{-x}\xi)\|_{\rm op}$. Also, $\OG(A(\zz))$ and $\OG_{H_\xi^\beta}(A(\zz))$ are the superoperators defined in Section~\ref{LipZd.ssect-supop}. Since $|h|^\beta\leq (1+h^2)^{\frac{\beta}{2}}$, the previous estimate reads as follows
\begin{align*}
\big|T_1(\zz)\varphi(x)\big|
	\leq & \frac{2^{1-\beta}\, C_L^\beta}{(r-\tR_H)^\beta} \, 
		\left( \big(
			\OG_{H_\xi^\beta}(A(\zz))\varphi
		\big)(x)+
		\big(
			H_\xi^\beta\, \OG(A(\zz))\varphi
		\big)(x)\right)\,.
\end{align*}

\noindent Lemma~\ref{LipZd.Lem-Hbeta} asserts that $H_\xi^\beta:\ell^2(\Ll)\to\ell^2(\Ll)$ is positivity preserving, self-adjoint satisfying $\|H_\xi^\beta\|\leq \|H\|_\beta$. Thus, Lemma~\ref{LipZd.Lem-G(A)} and Lemma~\ref{LipZd.Lem-G_D(A)} imply
\begin{align*}
\|T_1(\zz)\| 
	\;\leq\;  & \left(\frac{C_L}{2\, (r-\tR_H)}\right)^\beta \, 4 \, \ns \, \|H_\xi^\beta\| \, \|A(\zz)\|
	\;\leq\;  \left(\frac{C_L}{2\, (r-\tR_H)}\right)^\beta \, 4 \, \ns \, \|H\|_\beta \, \frac{1}{\dist\big(\zz,\sigma(H_\Theta)\big)}\,.
\end{align*}

\noindent Invoking the lower bound \eqref{eq:EsRes} on the distance $\dist\big(\zz,\sigma(H_\Theta)\big)$, the norm $\|T_1(\zz)\|$ is smaller or equal than $\frac{1}{4}$ since $C_{hop}\geq 1$, uniformly in $\xi\in\Xi$. 

\vspace{.1cm}

\noindent (iii): Let $\varphi\in\Hh$ and $x\in\Ll$. Estimate~\eqref{LipZd.Eq-EsRes} at the beginning of the proof leads to
\begin{align*}
\big|T_2(\zz)\varphi(x)\big|
	\;\leq\; &\sum_{z\in\Ll}\sum_{h\in\Rr}
		\left|\left(E_{h,z} \, U_h \, A_z(\zz) \, \Mch_z^{(r-\tR_H)} \, \varphi\right)(x)\right|\\
	\;\leq\; &\frac{2}{r^\beta} \, C_{hop}\, \sum_{z\in\Ll} \chi_z^{(r-\tR_H)}(x) 
		\sum_{h\in\Rr} \|t_h\|_\infty \, \Big|\big(U_h \, A_z(\zz) \, \Mch_z^{(r-\tR_H)} \, \varphi\big)(x)\Big|\,.
\end{align*}

\noindent As introduced in Lemma~\ref{LipZd.Lem-Hinfty}, $H^\infty:\ell^2(\Ll)\to\ell^2(\Ll)$ defined by $H^\infty:=\sum_{h\in\Rr} \|t_h\|_\infty \, U_h$ is a positivity preserving, self-adjoint operator. Hence, the previous considerations imply
$$\big|T_2(\zz)\varphi(x)\big|
	\;\leq \; \frac{2}{r^\beta} \, C_{hop} \, \big(\OG_{H^\infty}(A(\zz)) \, \varphi\big)(x)\,.
$$

\noindent Here $\OG_{H^\infty}(A(\zz)):\Hh\to\ell^2(\Ll)$ is again the superoperator introduced in Lemma~\ref{LipZd.Lem-G_D(A)}. Lemma~\ref{LipZd.Lem-Hinfty} additionally states that $\|H^\infty\|\leq \|H\|_\beta$. Hence, Lemma~\ref{LipZd.Lem-G_D(A)} leads to
$$\|T_2(\zz)\|
	\;\leq\; \frac{2}{r^\beta} \, C_{hop} \, \ns \, \|H^\infty\| \, \|A(\zz)\|
	\;\leq\; \frac{2}{r^\beta} \, C_{hop} \, \ns \, \|H\|_\beta \, \frac{1}{\dist\big(\zz,\sigma(H_\Theta)\big)}
	\;\leq\; \frac{1}{4}\,
$$

where the last inequality is a consequence of \eqref{eq:EsRes}. 

\vspace{.1cm}

\noindent (iv): Step (ii) and (iii) imply that $\|T_1(\zz)+T_2(\zz)\|\leq \frac{1}{2}$ uniformly in $\xi\in \Xi$, for all $\zz\in\rho(H_\Theta)$ satisfying \eqref{eq:EsRes}. Assuming for the moment that $\zz$ has a non-zero imaginary part, we know that $H_\xi-\zz$ is invertible because $H$ is self-adjoint and we can write:
$$
	(H_\xi-\zz)^{-1}=S(\zz)\left (\Id + T_1(\zz) + T_2(\zz)\right )^{-1}\,.
$$

The estimate \eqref{eq:Est-S} and $\|T_1(\zz)+T_2(\zz)\|\leq \frac{1}{2}$ imply: 
$$
\| (H_\xi-\zz)^{-1}\|
	\;\leq\; \frac{2\, \ns}{\dist\big(\zz,\sigma(H_\Theta)\big)}
$$

uniformly both in $\xi$ and in the imaginary part of $\zz$. By analytic continuation, the above estimate remains true for real elements of $\rho(H_\Theta)$ satisfying \eqref{eq:EsRes}. Moreover, the same estimate provides the uniform lower bound we are looking for:
$$\dist\big(\zz,\sigma(H_\xi)\big )
	\;\geq\; \frac{\dist\big(\zz,\sigma(H_\Theta)\big)}{2\, \ns}
	\,,\qquad \xi\in \Xi\,,
$$

which shows that $\zz\in \rho(H_\Xi)$.
\hfill $\Box$

\vspace{.2cm}
%

\noindent The previous proposition provides an estimate on the Hausdorff distance of the spectra whenever $\Hdco(\Xi,\Theta)$ is small enough. The following (classical) statement delivers an estimate if  $\Hdco(\Xi,\Theta)$ is greater or equal than a certain constant. A slightly more general version can be found in \cite[Chapter 5, Theorem 4.10]{Ka95}.

\begin{lemma}
\label{LipZd.Lem-p2Spec}
Let $A,B\in\Ll(\Hh)$ be self-adjoint. Then 
\begin{align*}
d_H\big(\sigma(A),\sigma(B)\big) \;
	\leq\; \|A-B\|
	\;\leq\; 2\, \max\big\{\|A\|\,,\,\|B\|\big\}\,.
\end{align*}
\end{lemma}

\noindent {\bf Proof: } Let $\lambda\not \in\sigma(A)$ such that $d(\lambda,\sigma(A))> \|A-B\|$. Then the operator $(B-A)(A-\lambda)^{-1}$ has norm strictly less than $1$ and so $\Id +(B-A)(A-\lambda)^{-1}$ is invertible. Thus
$$
	B-\lambda=\left (\Id +(B-A)(A-\lambda)^{-1}\right )(A-\lambda)
$$

is also invertible, which shows that $\lambda\not\in \sigma(B)$. In other words, no element of $\sigma(B)$ can be located at a distance larger than $\|A-B\|$ from $\sigma(A)$, which implies:
$$
\sup_{\lambda\in \sigma(B)} \dist\big(\lambda,\sigma(A)\big)
	\;\leq\; \|A-B\|\,.
$$

By interchanging $A$ with $B$, the proof is over. 
\hfill $\Box$

\subsection{Proof of Theorem~\ref{LipZd.th-main}} Recall the notation $\tR_H:= R_H+\|M^{-1}\|_{\max}\, \|M\|_{\max}+1$ and define the constant 
\begin{align*}
C_{d,\Ll}
	\;:=\; 16 \, \ns \, \max\big\{\|M^{-1}\|_{\max} \,\|M\|_{\max}\,,\, C_L \,,\, 1\big\} 
\end{align*}

\noindent which only depends on the lattice $\Ll$ and the choice of the Lipschitz-partition of unity (and hence on the dimension), c.f. Section~\ref{LipZd.ssect-qupa}. If $\Xi=\Theta$, then $\sigma(H_\Xi)=\sigma(H_\Theta)$. Now suppose $\Xi\neq\Theta$, namely $\Hdco(\Xi,\Theta)>0$. Set $r:=\Hdco(\Xi,\Theta)^{-1}$. We analyze two cases: (i) $1\leq r \leq 2\, \tR_H$ and (ii) $2\, \tR_H<r$. 

\vspace{.1cm}

\noindent (i): From $1\leq r \leq 2\, \tR_H$ we infer 
\begin{align*}
1 \;\leq\; \left(2\, \frac{\tR_H}{r}\right)^\beta
	\;\leq\; 6\, \max\{\|M^{-1}\|_{\max}\,\|M\|_{\max}\,,\,1\}\, \frac{R_H^\beta}{r^\beta}
\end{align*}

where the last inequality follows from the definition of $\tR_H$. Also, $H$ is self-adjoint and $\|H_{\xi'}\|\leq \|H\|_\beta$ for $\xi'\in\csp$. Then Lemma~\ref{LipZd.Lem-p2Spec} and the previous considerations imply
$$d_H\big(\sigma(H_\xi),\sigma(H_\eta)\big)
	\leq 2\, \|H\|_\beta
	\leq 12\, \max\{\|M^{-1}\|_{\max}\,\|M\|_{\max}\,,\,1\}\, R_H^\beta \|H\|_\beta \, \Hdco(\Xi,\Theta)^\beta
$$

\noindent for all $\xi\in\Xi$ and $\eta\in\Theta$. According to Remark~\ref{LipZd.Rem-ns3}, $\ns\geq 3$ always holds. Thus, the desired estimate on the Hausdorff distance of the spectra $\sigma(H_\Xi)$ and $\sigma(H_\Theta)$ is derived if $1\leq r \leq 2\, \tR_H$ as $C_{hop}\geq 1$.

\vspace{.1cm}

\noindent (ii): Suppose $r> 2\, \tR_H \geq 2$. For $\zz\in\sigma(H_\Xi)$, Proposition~\ref{LipZd.Prop-EsRes} leads to
$$
\dist\big(\zz,\sigma(H_\Theta)\big) \;
	\leq \; \frac{16 \, \ns \, \max\big\{C_L \,,\, 1\big\}}{2^\beta\, (r-\tR_H)^\beta}
		 \, C_{hop} \, \|H\|_\beta\,.
$$

\noindent By interchanging the role of $\Xi$ and $\Theta$ we obtain:
$$d_H\big(\sigma(H_\Xi),\sigma(H_\Theta)\big) \;
	\leq \; \frac{C_{d,\Ll}}{2^\beta} \, C_{hop} \, \|H\|_\beta \, \frac{1}{(r-\tR_H)^\beta}.
$$

\noindent In addition, the constraint $r>2 \, \tR_H$ implies $(r-\tR_H)^{-\beta}\leq \big(\frac{2}{r}\big)^\beta$ which together with $r:=\Hdco(\Xi,\Theta)^{-1}$ and the previous estimate finishes the proof.
\hfill$\Box$

\vspace{.2cm}

\subsection{Proof of Theorem~\ref{LipZd.th-infirange}}  If $\Xi=\Theta$, then $\sigma(H_\Xi)=\sigma(H_\Theta)$ follows. Thus, without loss of generality suppose $\Hdco(\Xi,\Theta)>0$. Let $r:=\Hdco(\Xi,\Theta)^{-1}\geq 1$. Recall that 
\begin{align*}
C_{d,\Ll}
	\;:=\; 16 \, \ns \, \max\big\{\|M^{-1}\|_{\max} \,\|M\|_{\max}\,,\, C_L \,,\, 1\big\} 
\end{align*}

Let us introduce a parameter $s\geq 1$ for which we consider the strongly pattern equivariant Hamiltonian $H|_s:=\big(H_{\xi'}|_s\big)_{\xi'\in\csp}$ with $\beta$-H\"older continuous coefficients, which is given by restricting the range to $\Rr\cap Q_s$, see Section~\ref{LipZd.ssect-ham}. 
Since the range of influence of $H$ has linear growth, there exists a $C_H\geq 1$ such that $R_{H|_s}$ (i.e. the radius of influence of $H|_s$) obeys $R_{H|_s}\leq C_H \, s$.

\vspace{.1cm}

First, assume that $1\leq r\leq 4 \, \|M^{-1}\|_{\max} \, \|M\|_{\max} + 4\, C_H$. The second inequality in Lemma~\ref{LipZd.Lem-p2Spec} implies
$$d_H\big(\sigma(H_\xi),\sigma(H_\eta)\big)
	\;\leq\; 2\, \|H\|_\beta
	\;\leq\; 2( 4 \, \|M^{-1}\|_{\max} \, \|M\|_{\max} + 4\, C_H)^\beta \, \|H\|_\beta \, \Hdco(\Xi,\Theta)^\beta
$$

\noindent for all $\xi\in\Xi$ and $\eta\in\Theta$. According to Remark~\ref{LipZd.Rem-ns3}, $\ns\geq 3$ always holds. This leads to the desired spectral estimate as $C_{hop}\geq 1$.

\vspace{.1cm}

\noindent Second, assume that $r> 4 \,\|M^{-1}\|_{\max} \, \|M\|_{\max} + 4\, C_H$. For every $s\geq 1$ we have
\begin{align*}
s^\beta \, \|H_{\xi'}-H_{\xi'}|_s\| 
	\leq \; s^\beta \sum_{|h|>s} \|t_h\|_\infty \leq \|H\|_\beta\, ,
\end{align*}

where we used that $s^\beta\leq (1+h^2)^{\frac{\beta}{2}}$ for $|h|>s$. Hence, $\|H_{\xi'}-H_{\xi'}|_s\|\leq s^{-\beta}\|H\|_\beta$ follows for all $\xi'\in\csp$. Since the operators are self-adjoint, Lemma~\ref{LipZd.Lem-p2Spec} implies 
\begin{align}\label{eq:infirange}
d_H\big(\sigma(H_{\xi'}),\sigma(H_{\xi'}|_s)\big) 
	\leq \|H_{\xi'}-H_{\xi'}|_s\|\leq \|H\|_\beta \, \frac{1}{s^{\beta}}
	\,, \qquad \xi'\in\csp\,.
\end{align}

\noindent Define $s:=\frac{r-2 \, \|M^{-1}\|_{\max} \, \|M\|_{\max}}{2 \, C_H} > 2$.  Then $R_{H|_s}\leq C_H \, s$ yields 
$$\frac{r}{2}
	\geq R_{H|_s} + \|M^{-1}\|_{\max} \, \|M\|_{\max}= \tR_{H|_s}-1\,.
$$

\noindent In particular, $r>\tR_{H|_s}$ holds. In addition, $r>4$ holds as $C_H\geq 1$ and so $r-\tR_{H|_s}>\frac{r}{4}$ follows by the previous considerations. Thus, Proposition~\ref{LipZd.Prop-EsRes} (applied to $H|_s$ and the given $r$) and $\|H|_s\|_\beta\leq \|H\|_\beta$ imply
$$d_H\big(\sigma(H_\Xi|_s),\sigma(H_\Theta|_s)\big) 
	\leq \frac{16\,\ns \max\{C_L\,,\, 1\}}{2^\beta \, \big(r-\tR_{H|_s}\big)^\beta} \, C_{hop} \, \|H|_s\|_\beta
	\leq \frac{2\, C_{d,\Ll}\, \ns }{r^\beta} \, C_{hop} \, \|H\|_\beta\,.
$$

\noindent The constraint $r>4 \,\|M^{-1}\|_{\max} \, \|M\|_{\max} + 4\, C_H$ leads to $\frac{r}{2}>2\,\|M^{-1}\|_{\max} \, \|M\|_{\max}$. Hence, from the definition of $s$, we obtain the estimate
$$\frac{1}{s}
	= \frac{2 \, C_H}{r-2\,\|M^{-1}\|_{\max}\, \|M\|_{\max}} 
	\leq  \frac{4 \, C_H}{r}.
$$

Since $\ns\geq 3$ holds in any dimension, we have $2\, 4^\beta\leq C_{d,\Ll}$. Combined with \eqref{eq:infirange}, the previous considerations imply
\begin{align*}
d_H\big(\sigma(H_\Xi),\sigma(H_\Theta)\big) \; 
	\leq \; &d_H\big(\sigma(H_\Xi),\sigma(H_\Xi|_s)\big) 
		+ d_H\big(\sigma(H_\Xi|_s),\sigma(H_\Theta|_s)\big) 
		+ d_H\big(\sigma(H_\Theta|_s),\sigma(H_\Theta)\big)\\
	\leq \; & 2 \,  \|H\|_\beta \, \frac{1}{s^{\beta}}
		+ 2\, C_{d,\Ll}\, \ns \, C_{hop} \, \|H\|_\beta \, \frac{1}{r^\beta}\\
	\leq \; &2\, C_{d,\Ll} \, (C_H^\beta + C_{hop}) \, \|H\|_\beta \; \Hdco\big(\Xi,\Theta\big)^\beta
\end{align*}

\noindent finishing the proof.
\hfill$\Box$


\end{document}